\documentclass[11pt,reqno]{amsart}
\usepackage[T1]{fontenc}
\usepackage{graphicx}
\usepackage{tikz}
\usepackage{cite}

\usepackage{color}
\definecolor{MyLinkColor}{rgb}{0,0,0.4}

\newcommand{\R}{{\mathbb R}}

\newcommand{\bA}{{\mathbb A}}
\newcommand{\bB}{\mathbb{B}}
\newcommand{\bT}{\mathbb{T}}
\newcommand{\Z}{{\mathbb Z}}
\newcommand{\N}{{\mathbb N}}

\newcommand{\kH}{\mathcal{H}}
\newcommand{\cO}{\mathcal{O}}

\newcommand{\kL}{\mathcal{L}}

\newcommand{\wh}{\widehat}
\newcommand{\wt}{\widetilde}

\newcommand{\re}{\mathop{\rm Re}\nolimits}
\newcommand{\im}{\mathop{\rm Im}\nolimits}
\newcommand{\PV}{\mathop{\rm PV}\nolimits}

\newcommand{\p}{\partial}

\newcommand{\e}{\varepsilon}

\newcommand{\0}{\Omega}

\newcommand{\supp}{\mathop{\rm supp}\nolimits}
\newtheorem{thm}{Theorem}[section]
\newtheorem{prop}[thm]{Proposition}
\newtheorem{lemma}[thm]{Lemma}

\theoremstyle{remark} 
\newtheorem{rem}[thm]{Remark}

\numberwithin{equation}{section} 

\title[The  nonlocal mean curvature flow of periodic graphs]{The  Nonlocal Mean Curvature Flow of Periodic Graphs}

 \author{Bogdan--Vasile Matioc}
\address{Fakult\"at f\"ur Mathematik, Universit\"at Regensburg,   93053 Regensburg, Deutschland.}
\email{bogdan.matioc@ur.de}

\author{Christoph Walker}
\address{Leibniz Universit\"at Hannover\\
Institut f\"ur Angewandte Mathematik\\
Welfengarten 1\\
30167 Hannover\\
Germany}
\email{walker@ifam.uni-hannover.de}


\usepackage[colorlinks=true,linkcolor=MyLinkColor,citecolor=MyLinkColor]{hyperref} 

\begin{document}

\begin{abstract}
We establish the well-posedness of the nonlocal mean curvature flow of order~${\alpha\in(0,1)}$  {for}  periodic graphs on $\R^n$ 
 {in all subcritical little H\"older spaces ${\rm h}^{1+\beta}(\bT^n)$ with $\beta\in(0,1)$}.
Furthermore, we prove that if the solution is initially sufficiently close to its integral mean in ${\rm h}^{1+\beta}(\bT^n)$, 
then it exists globally  in time  and converges exponentially fast towards a constant.
The proofs rely on the reformulation of the equation as a quasilinear evolution problem, which is shown to be of parabolic type by a direct localization approach, and on abstract parabolic theories for such problems.  
\end{abstract}
\maketitle

\section{Introduction and Main Results}\label {Sec:1}

Let $n\geq 1$ and  $ \Sigma:=\{(x,u(x))\,:\, x\in\R^n\}$ be the graph  of a function $u:\R^n\to\R$. 
The nonlocal fractional mean curvature $H_\alpha=H_\alpha(u)$ of order~${\alpha\in(0,1)}$  of $\Sigma$  is a quite recent notion introduced in \cite{CaffarelliSouganidis10}  and given by the formula
\[
H_\alpha(u)(x):=\frac{2}{\alpha}\int_{\Sigma} \frac{(y,u(y))-(x,u(x))}{|(y,u(y))-(x,u(x))|^{n+1+\alpha}}\cdot \nu(y)d  \sigma(y),\qquad x\in\R^n,
\]
 where 
\[
\nu:=(1+|\nabla u |^2)^{-1/2}(-\nabla u,1) \qquad\text{and}\qquad   d\sigma(y)=(1+|\nabla u |^2)^{1/2}(y)\,dy
\]
are the unit normal vector  and the volume element of $\Sigma$, respectively, and $a\cdot b$ denotes the Euclidean scalar product of two vectors $a$ and $b$. With the short hand notation
\[
\delta_{[x,y]}u:=u(x)-u(x-y),\qquad x,\, y\in\R^n,
\]
one may write the previous formula also as
\begin{align}\label{NMC}
H_\alpha(u)(x)=-\frac{2}{\alpha}\int_{\R^n} \frac{\delta_{[x,y]}u-y\cdot \nabla u(x-y)}{\big[|y|^2+(\delta_{[x,y]}u )^2\big]^{\tfrac{n+1+\alpha}{2}}}dy,\qquad x\in\R^n.
\end{align}
Besides this  {integral} representation of the nonlocal mean curvature there is also an alternative expression of $H_\alpha(u)$ as a principal value integral (see~\eqref{NMC2} below).
Actually, the nonlocal mean curvature may be defined more generally for (sufficiently smooth) open sets in~$\R^{n+1}$,  {see \cite{CaffarelliSouganidis10}},
 and the definition above corresponds to the particular case when the open set is  {a} subgraph. 
 In this general situation the nonlocal mean curvature has been a very active research topic in the recent past,  
 see for instance the review article  \cite{Fall18} and the references therein and in \cite{JulinLaManna20}. 

In this research, attention is focused on the nonlocal mean curvature flow
\begin{equation}\label{PB}
V(t)=-H_\alpha(u(t)),\quad t>0,\qquad \Sigma(0)=\Sigma_0\,,
\end{equation}
describing the evolution of the family of graphs $\{\Sigma(t)\,:\,t\ge 0\}$ with
\[
\Sigma(t)=\{(x,u(t,x))\,:\, x\in\R^n\},\quad t\ge 0,
\]
for $u(t,\cdot):\R^n\to\R$, where the initial graph $\Sigma_0:=\{(x,u_0(x))\,:\, x\in\R^n\}$ with $u_0:\R^n\to\R$ is given.
In \eqref{PB}, 
\begin{equation}\label{V}
V(t):=(1+|\nabla u(t)|^2)^{-1/2}\p_tu,\quad t>0,
\end{equation}
is the normal velocity of  $\Sigma(t)$. 

For the corresponding nonlocal mean curvature flow of sets, different concepts were used to study questions related to existence, uniqueness, regularity,  {occurence of neck pinching,
 and the singular limits when~${\alpha\to0}$  or $\alpha\to 1$}   {in the setting of viscosity} solutions, 
 see, e.g., \cite{CSV18, CP21, LKP22, Imbert09, CLNP21, CDNV19, CaffarelliSouganidis10, ChambolleMoriniPonsiglione15, Cameron19, CesaroniNovaga21}.
  Moreover, various properties  of smooth solutions (assuming their existence)  to the  fractional (volume preserving) mean curvature flow  
  were investigated in \cite{CSV20, SaezValdinoci19}.
 {Existence and   special features of hypersurfaces with constant and almost-constant nonlocal mean curvature have been established  in \cite{CFMN18, CFSW18, CFMW18a, CFMW18b, MNA20}.}

  The local well-posedness of the  {(volume preserving)} nonlocal mean curvature flow of  {bounded regular} sets 
  {in the setting of classical solutions} was  addressed  only recently in \cite{JulinLaManna20}  for initial data
 parametrizing~${\rm C}^{1,1}$-regular hypersurfaces.

Regarding the  nonlocal mean curvature flow \eqref{PB} of  {(possibly unbounded) graphs the local existence and uniqueness of classical
solutions determined by initial data with gradients in~${{\rm C}^\gamma(\R^n)}$, with~${\gamma>\alpha}$,}
  has recently been established in \cite{AFW22xx}  by means of analytic semigroup theory for quasilinear evolution equations. 
 {In both references \cite{JulinLaManna20, AFW22xx} the authors also showed that the solutions become instantaneously smooth.}

In this paper we consider the more specific situation of periodic graphs.
 Within this framework we reformulate~\eqref{PB} as a  quasilinear  {parabolic} evolution  {problem} in little H\"older spaces and show that it is well-posed and
  that the corresponding  {classical} solutions define a smooth semiflow in these spaces, see  {Theorem~\ref{MT1} below}.
   In the setting of classical solutions, this seems to be the first well-posedness result which covers all subcritical little H\"older spaces~${{\rm h}^{1+\beta}(\bT^n)}$ with $\beta\in(0,1)$, see Remark~\ref{R:1}.
   As in~\cite{AFW22xx} we rely on analytic semigroup theory,  {but} we  use a self-contained localization argument to derive the parabolicity of~\eqref{PB}.  
The semiflow property also enables us to  prove that the solution is smooth with respect to  both time and space variables for positive times and to investigate stability issues.
 Indeed, we shall prove in  {Theorem~\ref{MT2}} that if the datum~$u_0$ (defining the initial geometry $\Sigma_0$) is sufficiently close to its integral mean  {in the phase space},
  then the solution to~\eqref{PB} exists globally  in time and converges exponentially  {fast towards a constant (which possibly depends on $u_0$)}.

\subsection{Main Results}\label{Sec:1.2}

To be more precise, we first note from  {\eqref{NMC} and} \eqref{V} that the nonlocal mean curvature flow~\eqref{PB}  can be formulated as a quasilinear  {evolution}  problem
\begin{equation}\label{QAP}
\frac{du}{dt}(t)=\Phi(u(t))[u(t)], \quad t>0,\qquad u(0)=u_0,
\end{equation}
where 
\[
\Phi(u)[v](x):=\frac{2}{\alpha}(1+|\nabla u|^2)^{1/2}(x)\int_{\R^n} \frac{\delta_{[x,y]}v-y\cdot\nabla v(x-y)}{\big[|y|^2+(\delta_{[x,y]}u )^2\big]^{\tfrac{n+1+\alpha}{2}}}dy,\qquad x\in\R^n.
\]
Denoting by ${\rm h}^ {s}(\bT^n)$ the little H\"older  space on the torus  $\bT^n$ of order $ {s\geq0}$ (i.e. the closure of~${{\rm C}^ {\infty}(\bT^n)}$ in ${\rm C}^ {s}(\bT^n)$, see Section~\ref{Sec:1.1} below for a precise definition) and choosing the H\"older exponents~${\alpha,\,\beta,\,\gamma}$ such that 
\begin{equation}\label{constants}
 {\text{$\alpha,\, \beta\in(0,1)$ are arbitrary and  $\max\{\alpha,\, \beta\}<\gamma<\min\{1,\alpha+\beta\},$}}
\end{equation}
 we shall prove  that the evolution problem~\eqref{QAP} is of quasilinear parabolic type in the sense that
\begin{equation}\label{property}
-\Phi\in{\rm C}^\infty\big({\rm h}^{1+\beta}(\bT^n),\mathcal{H}({\rm h}^{1+\gamma}(\bT^n), {\rm h}^{ \gamma-\alpha}(\bT^n))\big),
\end{equation}
where  {the set}  $\mathcal{H}({\rm h}^{1+\gamma}(\bT^n), {\rm h}^{ \gamma-\alpha}(\bT^n))$ of   negative generators of analytic semigroups
   {is defined in Section~\ref{Sec:1.1}}. 
Property \eqref{property} enables us to use the quasilinear parabolic theory presented in \cite{Am93} (see also~\cite{MW20})  {in the context of \eqref{QAP}}.
 As a first main result we obtain the following theorem  ensuring the local well-posedness of~\eqref{QAP} 
and  characterizing the blow-up behavior of solutions which are not global. 
Moreover, it also states that the solution map defines a smooth semiflow on ${\rm h}^{1+\beta}(\bT^n)$ 
 {and that the solutions are smooth for positive times}.

\begin{thm}\label{MT1}
Let  {$\alpha,\,\beta,\,\gamma$ satisfy \eqref{constants}}.
Then, given~${u_0\in {\rm h}^{1+\beta}(\bT^n)}$, there exists a unique maximal solution $u:=u(\,\cdot\,; u_0)$ to \eqref{QAP} such that
\begin{equation*}  
 u\in {\rm C}([0,T^+),{\rm h}^{1+\beta}(\bT^n))\cap {\rm C}((0,T^+),{\rm h}^{1+\gamma}(\bT^n))\cap {\rm C}^1((0,T^+),{\rm h}^{\gamma-\alpha}(\bT^n)) ,
  \end{equation*}
with $T^+=T^+(u_0)\in(0,\infty]$ denoting  the maximal existence time.
 {Moreover,} the following properties hold:
\begin{itemize}
\item[(i)] If $T^+<\infty$, then
\[
\lim_{t\to T^+}\|u(t)\|_{1+\beta}=\infty.
\]
\item[(ii)] The map $[(t,u_0)\mapsto u(t;u_0)]$ defines a semiflow on ${\rm h}^{1+\beta}(\bT^n)$ which is smooth in the open set
  \[
\{(t,u_0)\,:\, u_0\in {\rm h}^{1+\beta}(\bT^n),\, 0<t<T^+(u_0)\}\subset \R\times {\rm h}^{1+\beta}(\bT^n)  
  \]
  and
    \begin{equation}\label{eq:fg}
 u\in {\rm C}^\infty((0,T^+)\times\bT^n).
  \end{equation}
  \item[(iii)] For $t\in[0,T^+)$ and $1\le j\le n$,
 \[
\|u(t)\|_0\leq \|u_0\|_0\qquad \text{and}\qquad
\|\partial_{x_j}u(t)\|_0\leq \|\partial_{x_j}u_0\|_0.
 \]
 {In the particular case $\beta>\alpha$ we  also have $\|\p_tu(t)\|_0\leq \|\p_t u(0)\|_0 $ for $t\in[0,T^+)$.}
\end{itemize}
 \end{thm}

 {We note} from (i) and (iii)  that the Hölder-seminorm of  {the gradient of a non-global solution blows  up in the limit $t\to T^+$.
Moreover, the additional assumption $\beta>\alpha$  imposed in (iii) for the estimate  $\|\p_tu(t)\|_0\leq \|\p_t u(0)\|_0 $ for $t\in[0,T^+)$ is needed in order 
to guarantee that $\p_t u(0)$ exists, as the integral \eqref{NMC} defining $H_\alpha(u)$ does in general not converge when $u\in {\rm C}^{1+\beta}(\bT^n)$ with $\beta\leq\alpha$.}

 {\begin{rem}\label{R:1}  We point out that if $u$ is a solution to \eqref{NMC}, then,  for $\lambda>0$, also the mapping~$u_\lambda$ given by
\[
u_\lambda(t,x):=\lambda^{-1}u(\lambda^{1+\alpha} t,\lambda x),
\]
is a solution to \eqref{NMC}. 
This scaling invariance property identifies ${\rm h}^1(\bT^n)$ as a critical space for the evolution problem~\eqref{NMC}. 
Hence, Theorem~\ref{MT1} covers all subcritical spaces ${{\rm h}^{1+\beta}(\bT^n)}$ with $\beta\in(0,1)$.
\end{rem}}

After having established in Section~\ref{Sec:2} suitable mapping properties of the operator $\Phi$ appearing in \eqref{QAP},
 we localize $\Phi(u)$ in Section~\ref{Sec:3}. 
 We then identify the obtained localized operators as Fourier multipliers and use the Mikhlin-H\"ormander theorem to establish their generator properties along with 
suitable resolvent estimates, 
see Proposition~\ref{A} . 
This enables us to establish the parabolicity of \eqref{QAP}, see Theorem~\ref{T:GP}, and implies thus \eqref{property}, so that Theorem~\ref{MT1} is  a consequence of the results from \cite{Am93}.
\\

In the second part of the paper we investigate stability properties of the flow. 
We show that the stationary solutions to~\eqref{QAP} are constant functions and derive the following stability result stating that any solution to~\eqref{QAP}, starting from an initial datum 
$u_0$ that is sufficiently close to its integral mean  {$\langle u_0\rangle$ in ${\rm h}^{1+\beta}(\bT^n)$},
 exists for all times and converges at an exponential rate  {towards   a constant} as $t\to\infty$.
 
   \begin{thm}\label{MT2}
   Let  {$\alpha,\,\beta,\,\gamma$ satisfy \eqref{constants}}.
Then, there is $\omega_0>0$ such that  for all~${\omega\in(0,\omega_0)}$  there are constants $\e>0$ and $M\geq1$ with the property that, for all $u_0\in {\rm h}^{1+\beta}(\bT^n)$ with
\[
\|u_0-\langle u_0\rangle\|_{1+\beta}\leq \e,
\]
the solution $u=u(\cdot; u_0)$ provided by Theorem~\ref{MT1} exists globally (that is $T^+(u_0)=+\infty$), and  there is a further constant $C(u_0)$ with $|C(u_0)|\leq \|u_0\|_0$ such that 
\begin{equation}\label{Exst}
\|u(t)-C(u_0)\|_{1+\beta}\leq Me^{-\omega t}\|u_0-\langle u_0\rangle\|_{1+\beta},\qquad t\geq0.
\end{equation}
   \end{thm}
A precise formula for  $C(u_0)$  can be found in the proof of Theorem~\ref{MT2}, see relation~\eqref{cu0}.
 Moreover, also the exponential rate $\omega_0$ is explicitly given in the proof of Theorem~\ref{MT2} (by~\eqref{repr} with $a=0$).

The proof of Theorem~\ref{MT2} is contained in Section~\ref{Sec:5}  and  relies on the principle of linearized stability of quasilinear evolution problems in interpolation spaces~\cite[Theorem~1.3]{MW20}. 
 {We emphasize that}  the nonlocal mean curvature flow~\eqref{PB} is invariant under  {vertical} translations of  graphs  {and therefore,}  the linearization of $\Phi$ 
at any stationary solution has  {zero as an eigenvalue}. 
 {To overcome this difficulty,} we  introduce  {in proof of Theorem~\ref{MT2}} a volume preserving unknown  {that solves a related} 
 quasilinear evolution problem,  {see \eqref{QAP'S},} to which the quasilinear principle of linearized stability established in \cite[Theorem~1.3]{MW20}  applies.
 {In this way we derive estimates  for the volume preserving unknown which are then used to estimate the actual solution $u=u(\cdot; u_0)$}.

\section{ {Notation} and mapping properties of $\Phi$} \label{Sec:2}

In this section we prove fundamental mapping properties of $\Phi$ in a suitable functional analytic setting  {which we now introduce}.

\subsection{Functional Analytic Setting}\label{Sec:1.1} 
Let  $E_1$ and $E_0$ be Banach spaces. We denote the space of bounded linear operators from $E_1$ to $E_0$ by $\kL(E_1,E_0)$ and set $\kL(E_0):=\kL(E_0,E_0)$. If $E_1$ is densely embedded in $E_0$,
 we~set (following~\cite{Am95})
\begin{equation*} 
\kH(E_1,E_0)=\{A\in\kL(E_1,E_0)\,:\, \text{$-A$ generates an analytic semigroup in $\kL(E_0)$}\}
\end{equation*}
 and note from~\cite[Theorem I.1.2.2]{Am95} that
\[
\kH(E_1,E_0)=\bigcup_{\substack{\kappa\geq1 \\ \omega>0}} \kH(E_1,E_0,\kappa,\omega),
\]
where $\kH(E_1,E_0,\kappa,\omega)$ consists of the operators $A$ with the property that  $\omega+A\in \kL(E_1,E_0)$ is an isomorphism and
\[
\kappa^{-1}\leq\frac{\|(\lambda+A)x\|}{|\lambda|\cdot\|x\|_{E_0}+\|x\|_{E_1}}\leq \kappa,\qquad 0\neq x\in E_1,\, \re\lambda\geq \omega.
\]

Given $k\in\N$, we  denote by ${\rm C}^k(\bT^n)$ the Banach space consisting of functions on $\R^n$ which are $2\pi$-periodic with respect to each variable and which posses continuous derivatives up to order  $k$.
Letting $\|\cdot\|_0$ be the supremum norm, the norm on   ${\rm C}^k(\bT^n)$  is as  usual defined by
\[
\|u\|_k  :=\max_{|\mu|\leq k}\|\p^\mu u\|_0.
\]
Given $\gamma\in(0,1),$ the H\"older space ${\rm C}^{k+\gamma}(\bT^n)$  consists of those functions $u\in{\rm C}^{k}(\bT^n)$ for which the norm
\[
\|u\|_{k+\gamma}:= [u ]_{k,\gamma} +\|u\|_k  
\]
is finite, where
\[
 [u ]_{k,\gamma} :=\max_{|\mu|=k}\sup_{x\neq y}\frac{|\p^\mu u(x)-\p^\mu u(y)|}{|x-y|^\gamma}.
\]
Finally, the little H\"older space ${\rm h}^ {k+\gamma}(\bT^n)$  is defined  as the closure of the set 
$${\rm C}^\infty(\bT^n):=\bigcap _{m\in\N}{\rm C}^m(\bT^n)$$ 
of smooth functions in
 ${\rm C}^{k+\gamma}(\bT^n)$ and ${\rm h}^ {k}(\bT^n):={\rm C}^{k}(\bT^n).$ 
  Then ${\rm h}^ {s}(\bT^n)$ is a Banach algebra for $s\ge0$.
If $s_2>s_1\geq0$, then ${\rm h}^ {s_2}(\bT^n)$ is a  densely and continuously embedded subspace of~${{\rm h}^ {s_1}(\bT^n)}$, since ${\rm C}^{s_2}(\bT^n)$ embeds continuously in ${\rm C}^{s_1}(\bT^n)$. Moreover, ${\rm C}^{s_2}(\bT^n)$ is a subspace of ${\rm h}^ {s_1}(\bT^n)$ whenever  $0\leq s_1<s_2$.

The little H\"older spaces are stable with respect to continuous interpolation; that is  for~${0\leq s_1<s_2}$ with $(1-\theta)s_1+\theta s_2\not\in\N$  we have
\begin{equation}\label{interpol}
({\rm h}^{s_1}(\bT^n), {\rm h}^{s_2}(\bT^n))_{\theta,\infty}^0={\rm h}^ {(1-\theta)s_1+\theta s_2}(\bT^n),
\end{equation}
where $(\cdot ,\cdot )_{\theta ,\infty}^0$ is the continuous interpolation functor of exponent $\theta\in(0,1),$ see \cite{DaPL88}.

\subsection{Mapping Properties of $\Phi$}

The goal of this section is to establish Proposition~\ref{P1} below.
  The mapping property  provided by this proposition together with the generator result for~$\Phi(u)$ established later in Theorem~\ref{T:GP} are the main steps in the proof of the aforementioned claim~\eqref{property}. 

\begin{prop}\label{P1}
If $\alpha,\,\beta,\,\gamma$  satisfy \eqref{constants}, then
\begin{equation}\label{MP}
\Phi\in{\rm C}^\infty\big({\rm h}^{1+\beta}(\bT^n),\mathcal{L}({\rm h}^{1+\gamma}(\bT^n), {\rm h}^{ \gamma-\alpha}(\bT^n))\big).
\end{equation}
\end{prop}

The proof of Proposition~\ref{P1} requires some preliminary investigations and is thus postponed to the end of this section.

We first observe that 
\begin{equation}\label{DEC}
\Phi(u)=(1+|\nabla u|^2)^{1/2}A(u),
\end{equation}
where
\begin{equation}\label{DEC0}
A(u)[v](x):=\frac{2}{\alpha}\int_{\R^n} \frac{\delta_{[x,y]}v-y\cdot\nabla v(x-y)}{|y|^{n+1+\alpha}}\frac{1}{\big[1+(\delta_{[x,y]}u/|y|)^2\big]^{\tfrac{n+1+\alpha}{2}}}dy,\qquad x\in\R^n.
\end{equation}
In order to prove Proposition~\ref{P1},  it is useful to view the operator $A$ as an element of a more general class which we now introduce.
Given $p\in\N$ and a smooth function $F:\R^2\to\R$,
we   define the linear operator~${A_{p}^F(u_1,u_2)[w_1,\ldots, w_p, \cdot]}$ by setting
\begin{align*}
A_{p}^F&(u_1,u_2)[w_1,\ldots, w_p, v](x)\\[1ex]
&:=\int_{\R^n} \frac{\delta_{[x,y]}v-y\cdot\nabla v(x-y)}{|y|^{n+1+\alpha}}F\bigg(\frac{\delta_{[x,y]}u_1}{|y|}, \frac{\delta_{[x,y]}u_2}{|y|}\bigg) 
\bigg(\prod_{j=1}^p \frac{\delta_{[x,y]}w_i}{|y|}\bigg)dy
\end{align*}
for $x\in\R^n.$ 
  Note that $A(u)[v]=A_{0}^{F}(u,0)[v]$ with $F(\xi_1,\xi_2):= \tfrac{2}{\alpha}(1+\xi_1^2)^{-\tfrac{n+1+\alpha}{2}} $.

The mapping property established in Lemma~\ref{L:1} is the first essential step in the proof of~\eqref{MP}.

\begin{lemma}\label{L:1}
 Given $p\in\N$,   mappings~${u_1,\, u_{2},\, w_j\in  {\rm C}^{1+\beta}(\bT^n),\, 1\leq j\leq p,}$ and a smooth function  $F:\R^2\to\R$,
 there exists a positive constant~${C=C(  p,\,F,\,\|u_1\|_{1+\beta}, \|u_2\|_{1+\beta})}$ such that 
\begin{equation}\label{prop1}
 \|A_{p}^F(u_1,u_2)[w_1,\ldots, w_p, \cdot]\|_{\mathcal{L}({\rm C}^{1+\gamma}(\bT^n), {\rm C}^{ \gamma-\alpha}(\bT^n))}\leq C\prod_{j=1}^p\|w_j\|_{1+\beta}.
\end{equation} 
\end{lemma}
\begin{proof}
Given $v\in{\rm C}^{1+\gamma}(\bT^n),$ we set for $x,\, y\in\R^n$ with $y\neq 0$
\[
Q(x,y):= \frac{\delta_{[x,y]}v-y\cdot\nabla v(x-y)}{|y|^{n+1+\alpha}}F\bigg(\frac{\delta_{[x,y]}u_1}{|y|}, \frac{\delta_{[x,y]}u_2}{|y|}\bigg) 
\bigg(\prod_{j=1}^p \frac{\delta_{[x,y]} {w_j}}{|y|}\bigg).
\]
Using the mean value theorem, we  then have 
\begin{equation}\label{f1}
\begin{aligned}
|\delta_{[x,y]}v-y\cdot\nabla v(x-y)|&=\bigg|\int_0^1y\cdot (\nabla v(x-y+sy)-\nabla v(x-y))\, ds\bigg|\\[1ex]
& \leq C\min\{ \|v\|_1 |y|\,,\, [ v]_{1,\gamma} | y|^{1+\gamma}\}\\[1ex]
&\leq  C\| v\|_{1+\gamma} \min\{ |y|\,,\, |y|^{1+\gamma}\},
\end{aligned}
\end{equation}
hence
\[
| Q(x,y)|\leq C\|v\|_{1+\gamma}\bigg(\prod_{j=1}^p\|w_j \|_{1}\bigg)\bigg(\frac{1}{|y|^{n+\alpha}}{\bf 1}_{\{|y|\geq1\}}(y)+\frac{1}{|y|^{n+\alpha-\gamma}}{\bf 1}_{\{|y|<1\}}(y)\bigg),
\]
where $C=C(F, \|u_1\|_1,\|u_2\|_1),$ and therefore
\begin{align*}
\bigg|\int_{\R^n} Q(x,y)\,dy\bigg|&\leq C\|v\|_{1+\gamma}\bigg(\prod_{j=1}^p\|w_j\|_{1}\bigg) \bigg(\int_{\{|y|\geq1\}}\frac{1}{|y|^{n+\alpha}}\, dy+\int_{\{|y|<1\}}\frac{1}{|y|^{n+\alpha-\gamma}}\, dy\bigg),
\end{align*}
the latter integrals being both finite  since $0<\alpha<\gamma$.
Hence,  $A_{p}^F(u_1,u_2)[w_1,\ldots, w_p, v]$ is well-defined, $2\pi$-periodic with respect to each variable,  and 
\begin{equation}\label{Est1}
\|A_{p}^F(u_1,u_2)[w_1,\ldots, w_p, v]\|_0\leq C\|v\|_{1+\gamma}\bigg(\prod_{j=1}^p\|w_j\|_{1}\bigg) .
\end{equation}
We next prove that $A_{p}^F(u_1,u_2)[w_1,\ldots, w_p, v]\in{\rm C}^{\gamma-\alpha}(\bT^n).$ Let therefore $x_1, x_2\in\R^n$  and~${y\not=0}$ be given.
We may estimate
\begin{align*}
|Q(x_1,y)- Q(x_2,y)|&\leq C\bigg[\bigg(\prod_{j=1}^p\|w_j\|_{1}\bigg)J_1(y)+J_2(y)\bigg],
\end{align*}
where  $C=C(F, \|u_1 \|_1,   \|u_2 \|_1)$ and
\begin{align*}
J_1(y)&:=\frac{\big|\delta_{[x_1,y]}v-y\cdot\nabla v(x_1-y)-\big[\delta_{[x_2,y]}v-y\cdot\nabla v(x_2-y)\big]\big|}{|y|^{n+1+\alpha}},\\[1ex]
J_2(y)&:=\frac{\big|\delta_{[x_2,y]}v-y\cdot\nabla v(x_2-y)\big|}{|y|^{n+1+\alpha}}\bigg[  \bigg(\prod_{j=1}^p\|w_j \|_{1}\bigg)\sum_{i=1}^2
\bigg|\frac{\delta_{[x_1,y]}u_i}{|y|}-\frac{\delta_{[x_2,y]}u_i}{|y|}\bigg|\\[1ex]
&\hspace{5cm} +\sum_{j=1}^p\bigg(\prod_{\substack{\ell=1 \\ \ell\neq j}}^p\|w_\ell \|_{1}\bigg)\bigg|\frac{\delta_{[x_1,y]}w_j}{|y|}-\frac{\delta_{[x_2,y]}w_j}{|y|}\bigg|\bigg].
\end{align*}
 Similarly as in \eqref{f1}, we  deduce that
 \begin{equation}\label{Est2A}
\begin{aligned}
J_1(y)&\leq\frac{1}{|y|^{n+\alpha}}\int_0^1\big|\nabla v(x_1-y+sy)-\nabla v(x_1-y)-\big[\nabla v(x_2-y+sy)-\nabla v(x_2-y)\big]\big|\, ds\\[1ex]
&\leq C\|v\|_{1+\gamma}\bigg(\frac{|x_1-x_2|^\gamma}{|y|^{n+\alpha}}{\bf 1}_{\{|y|\geq |x_1-x_2|\}}(y)+\frac{1}{|y|^{n+\alpha-\gamma}}{\bf 1}_{\{|y|<|x_1-x_2|\}}(y)\bigg),
\end{aligned}
\end{equation}
and therefore
\begin{equation}\label{Est2}
\begin{aligned}
\int_{\R^n}J_1(y)\, dy&\leq C\|v\|_{1+\gamma}\bigg(\int_{\{|y|\geq |x_1-x_2|\}}\frac{|x_1-x_2|^\gamma}{|y|^{n+\alpha}}\, dy+\int_{\{|y|<|x_1-x_2|\}}\frac{1}{|y|^{n+\alpha-\gamma}}\, dy\bigg)  \\[1ex]
&\leq C\|v\|_{1+\gamma}|x_1-x_2|^{\gamma-\alpha}.
\end{aligned}
\end{equation}
Moreover, 
\begin{align*}
J_2(y)&\leq\frac{1}{|y|^{n+\alpha}}\bigg(\int_0^1|\nabla v(x_2-y+sy) -\nabla v(x_2-y)|\, ds\bigg)\\[1ex]
&\qquad\times\bigg[ \bigg(\prod_{j=1}^p\|w_j \|_{1}\bigg)
\sum_{i=1}^2\int_0^1|\nabla u_i(x_1-y+sy)-\nabla u_i(x_2-y+sy)|\, ds\\[1ex]
&\hspace{1.25cm}\ +\sum_{j=1}^p\bigg(\prod_{ \substack{\ell=1 \\ \ell\neq j}}^p\|w_\ell \|_{1}\bigg)\int_0^1|\nabla w_j(x_1-y+sy)-\nabla w_j(x_2-y+sy)|\, ds \bigg]
\end{align*}
and, since $\beta>\gamma-\alpha$,  we deduce from \eqref{f1} that
\begin{align*}
J_2(y)&\leq C\bigg(\prod_{j=1}^p\|w_j \|_{1+\beta}\bigg)  \|v\|_{1+\gamma} |x_1-x_2|^{\gamma-\alpha}
\bigg(\frac{1}{|y|^{n+\alpha}}{\bf 1}_{\{|y|\geq 1\}}(y)+\frac{1}{|y|^{n+\alpha-\gamma}}{\bf 1}_{\{|y|<1\}}(y)\bigg),
\end{align*}
with $C=C(  p,\,F,\, \|u_1\|_{1+\beta},\, \|u_2\|_{1+\beta}).$
 Thus we have
\begin{equation}\label{Est3}
\int_{\R^n}J_2(y)\, dy\leq C\bigg(\prod_{j=1}^p\|w_j \|_{1+\beta}\bigg)\|v\|_{1+\gamma}|x_1-x_2|^{\gamma-\alpha}.
\end{equation}
Gathering \eqref{Est1}-\eqref{Est3}, we conclude \eqref{prop1}.
\end{proof}

Given Banach spaces $X$ and $Y$, we denote by ${\rm C}^{1-}(X,Y)$   the space of locally Lipschitz maps from $X$ to $Y$, and ${\rm C}^{\infty}(X,Y)$ is its subspace consisting of the smooth maps from $X$ to~$Y$.
Moreover, $\kL^p_{\rm sym}(X,Y)$, with $p\in\N,$ stands for the space of $p$-linear, bounded,  {and} symmetric maps $A:X^p\to Y$.

In order to prove the mapping property \eqref{MP} we introduce a particular  subclass of integral operators of the class $A_p^F$ introduced above.
Namely, given $p \in\N$ and a smooth function~${f:\R\to\R}$,  we set
\begin{align}\label{Afp}
A_{p}^f(u)[w_1,\ldots, w_p,v](x)&:=\int_{\R^n} \frac{\delta_{[x,y]}v-y\cdot\nabla v(x-y)}{|y|^{n+1+\alpha}}f\bigg( 
\frac{\delta_{[x,y]}u}{|y|}\bigg) \bigg(\prod_{j=1}^p \frac{\delta_{[x,y]}w_j}{|y|}\bigg)dy
\end{align}
for $x\in\R^n.$
Note that 
\begin{equation}\label{f2}
A(u)[v]=A_{0}^{f}(u)[v]\quad \text{with}\quad f(s)=\frac{2}{\alpha}(1+s^2)^{-\tfrac{n+1+\alpha}{2}}.
\end{equation}
Moreover, we have~${A_p^f(u)=A_p^F(u,0)}$  for $F(\xi_1,\xi_2):=f(\xi_1)$.
Therefore, we  infer from Lemma~\ref{L:1} that
\[
[u\mapsto A_p^f(u)]:{\rm C}^{1+\beta}(\bT^n)\to\kL^p_{\rm sym}({\rm C}^{1+\beta}(\bT^n),\kL({\rm C}^{1+\gamma}(\bT^n), {\rm C}^{\gamma-\alpha}(\bT^n)))
\]
 {is a well-defined mapping.}
The next lemma shows that this mapping is also Fr\'echet differentiable and its Fr\' echet derivative is defined by an integral operator from the same subclass  (with $f'$ denoting the derivative of $f$). 

\begin{lemma}\label{L:2}
 Given $p\in\N$ and a smooth function  $f:\R\to\R$, the mapping 
 \[
[u\mapsto A_p^f(u)]:{\rm C}^{1+\beta}(\bT^n)\to\kL^p_{\rm sym}({\rm C}^{1+\beta}(\bT^n),\kL({\rm C}^{1+\gamma}(\bT^n), {\rm C}^{\gamma-\alpha}(\bT^n)))
\]
defined in \eqref{Afp} is Fr\'echet differentiable and 
\begin{equation}\label{derapf}
\p A_p^f(u)[\wt u][w_1,\ldots,w_p,v]=  A_{p+1}^{f'}(u)[\wt u,w_1,\ldots, w_p,v]
\end{equation}
for all  {$u,\, \wt u,\, w_j\in {\rm C}^{1+\beta}(\R^n)$, $1\leq j\leq p$, and $v\in{\rm C}^{1+\gamma}(\R^n)$}.
\end{lemma}
\begin{proof}
Given  {$u,\, \wt u,\, w_j\in {\rm C}^{1+\beta}(\R^n)$, $1\leq j\leq p$, with~${\|\wt u\|_{1+\beta}\leq 1},$ and  $v\in {\rm C}^{1+\gamma}(\bT^n)$}, we have  
\begin{align*}
 &\big(A_p^f(u+\wt u) [w_1,\ldots,w_p,v]-A_p^f(u) [w_1,\ldots,w_p,v]- A_{p+1}^{f'}(u)[\wt u,w_1,\ldots, w_p,v]\big)(x)\\[1ex]
 &=\int_{\R^n}  R(x,y)\frac{\delta_{[x,y]}v-y\cdot\nabla v(x-y)}{|y|^{n+1+\alpha}}
  \bigg(\prod_{j=1}^p \frac{\delta_{[x,y]}w_j}{|y|}\bigg)dy,
\end{align*}
where, given $x,\, y\in\R^n$ with $y\neq 0$,  the mean  value theorem, yields 
\begin{align*}
R(x,y)&:= f\bigg( \frac{\delta_{[x,y]}(u+\wt u)}{|y|}\bigg)-f\bigg( \frac{\delta_{[x,y]}u}{|y|}\bigg)-\frac{\delta_{[x,y]}\wt u}{|y|}f'\bigg( \frac{\delta_{[x,y]}u}{|y|} \bigg)\\[1ex]
&=\bigg( \frac{\delta_{[x,y]}\wt u}{|y|}\bigg)^2\int_0^1\int_0^1 s f''\bigg(\frac{\delta_{[x,y]}(u+\tau s\wt u)}{|y|}\bigg)\, d\tau\, ds. 
\end{align*}
 Introducing the smooth function $F:\R^2\to\R$  by
\[
F(\xi_1,\xi_2):=\int_0^1\int_0^1 s f'' (\xi_1+\tau s\xi_2 )\, d\tau\, ds,
\]
 we get
\begin{align*}
 &A_p^f(u+\wt u) [w_1,\ldots,w_p,v]-A_p^f(u) [w_1,\ldots,w_p,v]- A_{p+1}^{f'}(u)[\wt u,w_1,\ldots, w_p,v]\\[1ex]
 &=A_{p+2}^F(u,\wt u) [\wt u,\wt u, w_1,\ldots,w_p,v].
\end{align*}
 {Keeping $u$ fixed,} we now infer from Lemma~\ref{L:1} that there exists a positive constant~$C$ such that for all functions~${\wt  u\in {\rm C}^{1+\beta}(\bT^n)}$ with $\|\wt u\|_{1+\beta}\leq 1$
 {and all $w_j\in {\rm C}^{1+\beta}(\R^n)$, $1\leq j\leq p$,} we have
\begin{align*}
 &\|A_p^f(u+\wt u) [w_1,\ldots,w_p,v]-A_p^f(u) [w_1,\ldots,w_p,v]- A_{p+1}^{f'}(u)[\wt u,w_1,\ldots, w_p,v]\|_{\gamma-\alpha}\\[1ex]
 &\leq C\|v\|_{1+\gamma}\|\wt u\|_{1+\beta}^2\prod_{j=1}^p\|w_j\|_{1+\beta},
\end{align*}
and the claim follows.
\end{proof}

 We are now in a position to provide the proof of Proposition~\ref{P1}.

\begin{proof}[Proof of Proposition~\ref{P1}]
We first recall that $\Phi(u)=(1+|\nabla u|^2)^{1/2}A(u),$ see \eqref{DEC},
where, as a direct consequence of Lemma~\ref{L:2} and \eqref{f2}, we have 
\[
A\in{\rm C}^\infty({\rm C}^{1+\beta}(\bT^n),\mathcal{L}({\rm C}^{1+\gamma}(\bT^n), {\rm C}^{ \gamma-\alpha}(\bT^n))).
\]
If $u,\,v \in {\rm C}^{\infty}(\bT^n)$,  the latter property implies that  $A(u)[v]\in {\rm C}^{ \gamma-\alpha+\e}(\bT^n)$ for all $\e\in(0,1-\gamma)$.
But, as  ${\rm C}^{ \gamma-\alpha+\e}(\bT^n)$ is subspace of ${\rm h}^{ \gamma-\alpha}(\bT^n)$, we have $A(u)[v]\in {\rm h}^{ \gamma-\alpha}(\bT^n)$.
The smoothness of $A$ now ensures that 
\[
A\in{\rm C}^\infty({\rm h}^{1+\beta}(\bT^n),\mathcal{L}({\rm h}^{1+\gamma}(\bT^n), {\rm h}^{ \gamma-\alpha}(\bT^n))).
\]
Finally, observing that 
\[
\big[u\mapsto (1+|\nabla u|^2)^{1/2}\big]\in{\rm C}^\infty({\rm h}^{1+\beta}(\bT^n), {\rm h}^{\beta}(\bT^n))
\]
and taking into account that $\beta>\gamma-\alpha$, we obtain the desired regularity property for~$\Phi$.
\end{proof}

\section{The Parabolicity Property} \label{Sec:3}
 {Throughout this section we again assume that the H\"older exponents $\alpha,\,\beta,\, \gamma$ satisfy \eqref{constants}.} 
The main goal is to show that, given $u\in{\rm h}^{1+\beta}(\bT^n)$, the operator $\Phi(u)$, viewed as an unbounded operator in $ {\rm h}^{\gamma-\alpha}(\bT^n)$ 
with {domain of definition} $ {\rm h}^{1+\gamma}(\bT^n)$, 
generates an analytic semigroup  in $\kL({\rm h}^{\gamma-\alpha}(\bT^n))$. 
 This is the content of the next theorem.
\begin{thm}\label{T:GP}
If $u\in{\rm h}^{1+\beta}(\bT^n) $, then $-\Phi(u)\in\mathcal{H}({\rm h}^{1+\gamma}(\bT^n), {\rm h}^{ \gamma-\alpha}(\bT^n)).$
\end{thm}
The proof of Theorem~\ref{T:GP} requires some preparation and is therefore postponed to the end of this section.
 As a starting point we define the continuous path
\[
[\tau\mapsto \Phi(\tau u)]:[0,1]\to \mathcal{L}({\rm h}^{1+\gamma}(\bT^n), {\rm h}^{ \gamma-\alpha}(\bT^n))
\]
connecting $\Phi(u)$ to the operator $\Phi(0)$.
In order to establish Theorem~\ref{T:GP}, we follow a strategy inspired  by \cite{E94,  ES95,  ES97}. 
 In a first step we approximate locally the operator $\Phi(\tau u)$ by certain operators with ``freezed'' $u$, see Proposition~\ref{P:1} below. 
As a second step we identify these operators as Fourier multipliers and use the Mikhlin-H\"ormander theorem to establish their generator properties along with 
 (uniform) resolvent estimates, see \eqref{uniest} in Proposition~\ref{A}. 
 Together with 
Proposition~\ref{P:1} this builds the core of the proof of  Theorem~\ref{T:GP}.

\subsection{Localization of $\Phi(\tau u)$} \label{Sec:3.1}
Guided by~\eqref{DEC0}, we  start by introducing  the class of operators  
\begin{equation}\label{fourmult}
A^a[v](x):=\frac{2}{\alpha}\int_{\R^n} \frac{\delta_{[x,y]}v-y\cdot\nabla v(x-y)}{|y|^{n+1+\alpha}}\frac{1}{\big[1+\tau^2\big(| y\cdot a|/|y|)^2\big]^{\tfrac{n+1+\alpha}{2}}}dy
\end{equation}
 for $a\in\R^n$, $v\in {\rm h}^{1+\gamma}(\bT^n)$,  and~${x\in\R}$.
It is straightforward to infer from the arguments used in the first part of the proof of Lemma~\ref{L:1} that~${A^a\in \mathcal{L}({\rm h}^{1+\gamma}(\bT^n), {\rm h}^{ \gamma-\alpha}(\bT^n))}$ with
\begin{equation}\label{estAa}
\sup_{a\in\R^n}\|A^a\|_{\mathcal{L}({\rm h}^{1+\gamma}(\bT^n), {\rm h}^{ \gamma-\alpha}(\bT^n))}<\infty.
\end{equation}
We shall see later on that $A^a$  is  a Fourier multiplier and that its symbol  can be expressed explicitly in terms of an integral, see \eqref{fmult}-\eqref{symbol}  {below}.
It is worthwhile  to point out the relation
\begin{equation}\label{phi0}
A^0=\Phi(0).
\end{equation}

In order to locally approximate the operator $\Phi(\tau u)$, $\tau\in[0,1]$, by Fourier multipliers $\delta A ^a$ (with $\delta\geq1$ and $a\in\R^n$),
we choose  for each~${\varepsilon\in(0,1)}$ a finite $\varepsilon$-localization family; that is,  a family  $\{(\pi_j^\varepsilon, x_j^\e)\,:\, 0\leq j\leq N\}\subset  {\rm C}^\infty(\bT^n,[0,1])\times\R,$
with $N=N(\varepsilon)\in\mathbb{N} $ sufficiently large such that 
\begin{align*}
\bullet\,\,\,\, \,\,  & \mbox{$\supp \pi_j^\varepsilon =\overline{\bB}_\e(x_j^\e)+2\pi\Z^n$ for $0\leq j\leq N$,} \\[1ex]
\bullet\,\,\,\, \,\, &\mbox{ $\displaystyle\sum_{j=0}^N \pi_j^\varepsilon=1$ on $\bT^n$.} 
\end{align*} 
 With such a finite $\varepsilon$-localization family we associate  a family   
$$
\{\chi_j^\varepsilon\,:\, 0\leq j\leq N\}\subset {\rm C}^\infty(\bT^n,[0,1])
$$  
satisfying
\begin{align*}
\bullet\,\,\,\, \,\,  &\mbox{$\chi_j^\varepsilon=1$ on $\supp \pi_j^\varepsilon$,} \\[1ex]
\bullet\,\,\,\, \,\,  &\mbox{$\supp \chi_j^\varepsilon \subset \overline{\bB}_{3\e}(x_j^\e)+2\pi\Z^n$ for $0\leq j\leq N$.} 
\end{align*} 
It readily follows  {from the above properties} that, for each $r\ge 0$, the map 
\begin{align}\label{eqnorm}
\bigg[f\mapsto \sum_{j=0}^N\|\pi_j^\varepsilon f\|_{r}\bigg]: {\rm h}^r(\bT^n)\to[0,\infty)
\end{align}
defines an equivalent norm on ${\rm h}^r(\bT^n)$.
We are now in a position to present the aforementioned localization result.

 \begin{prop}\label{P:1}
Let  {$\gamma'\in(\max\{\alpha,\beta\},\gamma)$} and let   $u\in {\rm h}^{1+\beta}(\bT^n)$ and~$\nu>0$ be given. 
Then, if $\e\in(0,1)$ is sufficiently small, there exists a positive constant~${K=K(\e)}$ such that for all~${\tau\in[0,1]},$ $v\in {\rm h}^{1+\gamma}(\bT^n)$, and $0\leq j\leq N$ we have
 \begin{equation}\label{Dj}
  \big\|\pi_j^\e \Phi(\tau u)[v]- (1+\tau^2|\nabla u|^2)^{1/2}(x_j^\e)A^{\tau \nabla u(x_j^\e)}[\pi_j^\e v]\big\|_{\gamma-\alpha}\leq \nu \|\pi_j^\e v\|_{1+\gamma}+K\|  v\|_{1+\gamma'}.
 \end{equation}
\end{prop}

The proof of Proposition~\ref{P:1} is postponed to the end of this subsection as it is based  on the auxiliary results   presented in Lemma~\ref{L:4P1} and Lemma~\ref{L:5P1}  below.
We start with an estimate for the commutator $[\pi_j^\e,\Phi(\tau u)].$

 \begin{lemma}\label{L:4P1}
 {Let}   $u\in {\rm h}^{1+\beta}(\bT^n)$ and~${\e\in(0,1)}$ be given. 
Then,  there is a positive constant~${K=K(\e)}$ such that for all~${\tau\in[0,1],\, v\in {\rm h}^{1+\gamma}(\bT^n)}$, and $0\leq j\leq N$ we have
 \begin{equation}\label{DjA1}
\|\pi_j^\e \Phi(\tau u)[v]- \Phi(\tau u)[\pi_j^\e v]\|_{\gamma-\alpha}\leq K\|  v\|_{1+\gamma-\alpha}.
 \end{equation}
\end{lemma}
\begin{proof}
In this proof constants denoted by $C$ depend only on $u$  and $\alpha,\,\beta,\,\gamma$, while constants denoted by $K$ may depend additionally also on $\e$. 

Letting  $\varphi:=\pi_j^\e \Phi(\tau u)[v]- \Phi(\tau u)[\pi_j^\e v],$ we have
\[
\varphi=-\frac{2}{\alpha}(1+\tau^2|\nabla u|^2)^{1/2}(\varphi_1+\varphi_2),
\]
where, for $x\in\R^n$, 
\begin{align*}
\varphi_1(x)&:=\int_{\R^n} \frac{\delta_{[x,y]} \pi_j^\e  }{|y|^{n+1+\alpha}}\frac{y\cdot\nabla  v (x-y)}{\big[1+\tau ^2(\delta_{[x,y]}u/|y|)^2\big]^{\tfrac{n+1+\alpha}{2}}}dy,\\[1ex]
\varphi_2(x)&:=\int_{\R^n} \frac{\delta_{[x,y]} \pi_j^\e-y\cdot\nabla \pi_j^\e(x-y)  }{|y|^{n+1+\alpha}}\frac{  v (x-y)}{\big[1+\tau ^2(\delta_{[x,y]}u/|y|)^2\big]^{\tfrac{n+1+\alpha}{2}}}dy.
\end{align*}
 Since $\beta> \gamma-\alpha$ and and the little H\"older spaces are Banach algebras,
 we have 
 \begin{equation}\label{bvb1}
 \|\varphi\|_{\gamma-\alpha}\leq C\big(\|\varphi_1\|_{\gamma-\alpha}+\|\varphi_2\|_{\gamma-\alpha}\big).
 \end{equation}
We next estimate $\|\varphi_1\|_{\gamma-\alpha}$.
Using the mean value theorem,   for $x,\, y\in\R^n$ with~${y\neq 0}$ we have
 \[
 \Bigg| \frac{\delta_{[x,y]} \pi_j^\e  }{|y|^{n+1+\alpha}}\frac{y\cdot\nabla  v (x-y)}{\big[1+\tau ^2(\delta_{[x,y]}u/|y|)^2\big]^{\tfrac{n+1+\alpha}{2}}}\Bigg|
 \leq  {C}\frac{\|v\|_1}{|y|^{n+\alpha}}\big({\bf 1}_{\{|y|\geq1\}}(y)+\|\pi_j^\e\|_1|y| {\bf 1}_{\{|y|<1\}}(y)\big),
 \]
 and therefore
  \begin{equation}\label{bvb2}
 \|\varphi_1\|_{0}\leq K\| v\|_{1}.
 \end{equation}
 In order to estimate the seminorm $[\varphi_1]_{\gamma-\alpha},$ let $x_1,\, x_2\in\R^n$ be given. Then, the mean value theorem yields
\begin{align*}
|\varphi_1(x_1)-\varphi_1(x_2)|&\leq\int_{\R^n}\frac{\big|(\delta_{[x_1,y]} \pi_j^\e)\nabla v (x_1-y) -(\delta_{[x_2,y]} \pi_j^\e)\nabla v(x_2-y)\big| }{|y|^{n+\alpha}}dy\\[1ex]
&\quad+C\int_\R \frac{\big|(\delta_{[x_2,y]} \pi_j^\e)\nabla v (x_2-y)\big| }{|y|^{n+\alpha}}\Big| \frac{\delta_{[x_1,y]}u}{y}-\frac{\delta_{[x_2,y]}u}{y} \Big|dy.
\end{align*}
Taking again advantage of the mean value theorem, we obtain
\begin{align*}
&\big|(\delta_{[x_1,y]} \pi_j^\e)\nabla v (x_1-y) -(\delta_{[x_2,y]} \pi_j^\e)\nabla v(x_2-y)\big|\\[1ex]
&\leq C\|v\|_1\big| \delta_{[x_1,y]} \pi_j^\e  - \delta_{[x_2,y]} \pi_j^\e \big|+\big|(\delta_{[x_2,y]}\pi_j^\e)(\nabla v (x_1-y) -\nabla v(x_2-y))\big| \\[1ex]
&\leq C\|v\|_1\big([\pi_j^\e]_{1,\gamma-\alpha}|y|{\bf 1}_{\{|y|<1\}}(y)+[\pi_j^\e]_{\gamma-\alpha}{\bf 1}_{\{|y|\geq1\}}(y)\big)|x_1-x_2|^{\gamma-\alpha}\\[1ex]
&\quad+ C\|v \|_{1+\gamma-\alpha}\big(\|\pi_j^\e\|_1|y|{\bf 1}_{\{|y|<1\}}(y)+{\bf 1}_{\{|y|\geq1\}}(y)\big)|x_1-x_2|^{\gamma-\alpha}
\end{align*}
and, in view of $\beta>\gamma-\alpha$, 
\begin{align*}
&\big|(\delta_{[x_2,y]} \pi_j^\e)\nabla v (x_2-y)\big(\delta_{[x_2,y]}u-\delta_{[x_1,y]}u\big) \big| \\[1ex]
&\leq C\|v\|_1 [u ]_{1,\gamma-\alpha}\big(\|\pi_j^\e\|_{1}|y|^2{\bf 1}_{\{|y|<1\}}(y)+|y|{\bf 1}_{\{|y|\geq1\}}(y)\big)|x_1-x_2|^{\gamma-\alpha}.
\end{align*}
The latter estimates lead us to
 \begin{equation}\label{bvb3}
 [\varphi_1]_{\gamma-\alpha}\leq K\|v\|_{1+\gamma-\alpha}.
 \end{equation}
Gathering \eqref{bvb2}-\eqref{bvb3}, we get
 \begin{equation}\label{bvb4}
 \|\varphi_1\|_{\gamma-\alpha}\leq K\|v\|_{1+\gamma-\alpha}.
 \end{equation}
 
 Concerning $\varphi_2$, we note  for $x\in\R^n$ that
\begin{align*}
\varphi_2(x)=\int_0^1\int_{\R^n} \frac{y\cdot(\nabla \pi_j^\e(x-y+sy)-\nabla \pi_j^\e(x-y))}{|y|^{n+1+\alpha}}\frac{  v (x-y)}{\big[1+\tau ^2(\delta_{[x,y]}u/|y|)^2\big]^{\tfrac{n+1+\alpha}{2}}}\, dy\,ds,
\end{align*}
and we  therefore  may repeat the arguments used to derive \eqref{bvb4} and conclude that   
\begin{equation}\label{bvb5}
 \|\varphi_2\|_{\gamma-\alpha}\leq K\|v\|_{1+\gamma-\alpha}.
 \end{equation}
The desired claim follows now from \eqref{bvb1}, \eqref{bvb4}, and \eqref{bvb5}.
\end{proof}

Lemma~\ref{L:5P1}  is the second main ingredient in the proof of Proposition~\ref{P:1}:

\begin{lemma}\label{L:5P1}
Let  {$\gamma'\in(\max\{\alpha,\beta\},\gamma)$} and let   $u\in {\rm h}^{1+\beta}(\bT^n)$ and~$\nu>0$ be given. 
Then, if $\e\in(0,1)$ is sufficiently small, there exists a positive constant~${K=K(\e)}$ such that for all~${\tau\in[0,1]},$ $v\in {\rm h}^{1+\gamma}(\bT^n)$, and $0\leq j\leq N$ we have
 \begin{equation}\label{Dj2}
  \big\| \Phi(\tau u)[\pi_j^\e v]- (1+\tau^2|\nabla u|^2)^{1/2}(x_j^\e)A^{\tau \nabla u(x_j^\e)}[\pi_j^\e v]\big\|_{\gamma-\alpha}\leq \nu \|\pi_j^\e v\|_{1+\gamma}+K\|  v\|_{1+\gamma'}.
 \end{equation}
\end{lemma}
\begin{proof}
As before, constants denoted by $C$ depend only on $u$  and $\alpha,\,\beta,\,\gamma,\,  {\gamma'}$,  while constants denoted by $K$ may depend additionally also on $\e$. 
Let
\[
\psi:=\Phi(\tau u)[\pi_j^\e v]- (1+\tau^2|\nabla u|^2)^{1/2}(x_j^\e)A^{\tau \nabla u(x_j^\e)}[\pi_j^\e v].
\]
Taking into account that $\chi_j^\e\pi_j^\e=\pi_j^\e$ and recalling~\eqref{DEC}, we  have
\[
\psi=\psi_1-\psi_2,
\]
where
\begin{align*}
\psi_1&:=  \chi_j^\e \big((1+\tau^2 |\nabla u|^2)^{1/2}A(\tau u)[\pi_j^\e v]-  (1+\tau^2|\nabla u|^2)^{1/2}(x_j^\e)A^{\tau \nabla u(x_j^\e)}[\pi_j^\e v]\big),\\[1ex]
\psi_2&:=  [\chi_j^\e,\Phi(\tau u)][\pi_j^\e v] -(1+\tau^2|\nabla u|^2)^{1/2}(x_j^\e)[\chi_j^\e, A^{\tau \nabla u(x_j^\e)}][\pi_j^\e v].
\end{align*}
It is not difficult to see that the arguments in the proof of Lemma~\ref{L:4P1} can also  be used to estimate  the commutators $[\chi_j^\e,\Phi(\tau u)]$ and $[\chi_j^\e, A^{\tau \nabla u(x_j^\e)}]$
which appear in the definition of~$\psi_2$ above. Hence
 \begin{equation}\label{96a}
 \|\psi_2\|_{\gamma-\alpha}\leq K\|\pi_j^\e v\|_{1+\gamma-\alpha}\leq K\|v\|_{1+\gamma'}
 \end{equation}
 as $\gamma'>\beta>\gamma-\alpha$.
It remains to estimate the function
\begin{align*}
\psi_1(x)=\chi_j^\e(x)\int_{\R^n}  P(x,y)dy,\qquad x\in\R^n,
\end{align*}
where, given $x,\, y\in\R^n$ with $y\neq0$, we set
\begin{align*}
P(x,y)&:=\frac{2}{\alpha}\frac{\delta_{[x,y]}(\pi_j^\e v)-y\cdot\nabla (\pi_j^\e v)(x-y)}{|y|^{n+1+\alpha}}\\[1ex]
&\qquad\times\Bigg(\frac{(1+\tau^2 |\nabla u|^2)^{1/2}(x)}{\big[1+\tau^2(\delta_{[x,y]}u/|y|)^2\big]^{\tfrac{n+1+\alpha}{2}}}
-\frac{(1+\tau^2 |\nabla u|^2)^{1/2}(x_j^\e)}{\big[1+\tau^2\big(|y\cdot \nabla u(x_j^\e)|/|y|)^2\big]^{\tfrac{n+1+\alpha}{2}}}\Bigg).
\end{align*}

We first infer from  \eqref{f1} (with $\gamma$ replaced by $\gamma'$) that
\begin{equation}\label{96b}
 \|\psi_1\|_0\leq C\|\pi_j^\e v\|_{1+\gamma'}\leq K\|v\|_{1+\gamma'}.
 \end{equation}
Let now $x_1\neq x_2\in\R^n$. 
In order to estimate the difference $|\psi_1(x_1)-\psi_1(x_2)| $ we may assume, 
in view of the property $\supp \chi_j^\varepsilon \subset \overline{\bB}_{3\e}(x_j^\e)+2\pi\Z^n$, that $x_1,\, x_2\in \overline{\bB}_{3\e}(x_j^\e)$.
We then have
\[
|\psi(x_1)-\psi(x_2)| \leq S_1+S_2,
\]
where
\begin{align*}
 S_1:=\int_{\R^n} | P(x_1,y)- P(x_2,y)|\, dy\qquad\text{and}\qquad S_2:=|\chi_j^\e(x_1)-\chi_j^\e(x_2)|\int_{\R^n} |P(x_2,y)|\, dy.
\end{align*}
 The estimate \eqref{f1} (with $\gamma$ replaced by $\gamma'$)  yields
 \begin{equation}\label{96c}
 S_2\leq C[\chi_j^\e]_{\gamma-\alpha} \|\pi_j^\e v\|_{1+\gamma'}|x_1-x_2|^{\gamma-\alpha}\leq K\|v\|_{1+\gamma'}|x_1-x_2|^{\gamma-\alpha},
 \end{equation}
 and it remains to estimate the term $ S_1$.
To this end we use the mean value theorem to derive
\[
| P(x_1,y)- P(x_2,y)|\leq C(L_1(y)+ L_2(y) {+L_3(y)}),\qquad y\in\R^n\setminus\{0\},
\]
where
\begin{align*}
 L_1(y)&:= \frac{\big|\delta_{[x_1,y]}(\pi_j^\e v)-y\cdot\nabla (\pi_j^\e v)(x_1-y)-[\delta_{[x_2,y]}(\pi_j^\e v)-y\cdot\nabla (\pi_j^\e v)(x_2-y)]\big|}{|y|^{n+1+\alpha}}\\[1ex]
&\qquad\times|\nabla u(x_1)-\nabla u(x_j^\e)|,\\[1ex]
L_2(y)&:= \frac{\big|\delta_{[x_1,y]}(\pi_j^\e v)-y\cdot\nabla (\pi_j^\e v)(x_1-y)-[\delta_{[x_2,y]}(\pi_j^\e v)-y\cdot\nabla (\pi_j^\e v)(x_2-y)]\big|}{|y|^{n+1+\alpha}}\\[1ex]
&\qquad\times\bigg|\frac{\delta_{[x_1,y]}u}{|y|}-\frac{y\cdot \nabla u(x_j^\e)}{|y|}\bigg|,\\[1ex]
L_3(y)&:= \frac{\big|\delta_{[x_2,y]}(\pi_j^\e v)-y\cdot\nabla (\pi_j^\e v)(x_2-y) \big|}{|y|^{n+1+\alpha}}
\Bigg(|\nabla u(x_1)-\nabla u(x_2)|+\bigg|\frac{\delta_{[x_1,y]}u}{|y|}-\frac{\delta_{[x_2,y]}u}{|y|}\bigg|\Bigg).
\end{align*}
Recalling \eqref{Est2A},   we have in view of $u\in{ \rm h}^{1+\beta}(\bT^n) $ and $x_1\in \overline{\bB}_{3\e}(x_j^\e)$  
\begin{equation*}
L_{1}(y)\leq C\e^\beta\|\pi_j^\e v\|_{1+\gamma}\bigg(\frac{|x_1-x_2|^\gamma}{|y|^{n+\alpha}}{\bf 1}_{\{|y|\geq |x_1-x_2|\}}(y)+\frac{1}{|y|^{n+\alpha-\gamma}}{\bf 1}_{\{|y|<|x_1-x_2|\}}(y)\bigg).
\end{equation*}
Hence, if $\e$ is sufficiently small, then
\begin{equation}\label{96d}
 \int_{\R^n}C  L_1(y)\, dy \leq \frac{\nu}{2} \|\pi_j^\e v\|_{1+\gamma}|x_1-x_2|^{\gamma-\alpha}
 \end{equation}
 for all~${\tau\in[0,1]},$ $v\in {\rm h}^{1+\gamma}(\bT^n)$, and $0\leq j\leq N$.
Moreover, since $|x_1-x_2| \leq 6\e$,   the estimate~\eqref{Est2A}  and the mean value theorem  lead us to
\begin{align*}
L_2(y)&\leq C\e^\beta\|\pi_j^\e v\|_{1+\gamma}\bigg(\frac{|x_1-x_2|^\gamma}{|y|^{n+\alpha}}{\bf 1}_{\{|x_1-x_2|\leq |y|\}}(y)+\frac{1}{|y|^{n+\alpha-\gamma}}{\bf 1}_{\{|y|<|x_1-x_2|\}}(y)\bigg)
{\bf 1}_{\{|y| {\leq}6\e\}}(y)\\[1ex]
&\quad+C\frac{1}{|y|^{n+\alpha}}[\pi_j^\e v]_{1,\gamma-\alpha}|x_1-x_2|^{\gamma-\alpha}{\bf 1}_{\{|y| {>} 6\e\}}(y)\\[1ex]
&\leq C\e^\beta\|\pi_j^\e v\|_{1+\gamma}\bigg(\frac{|x_1-x_2|^\gamma}{|y|^{n+\alpha}}{\bf 1}_{\{ |x_1-x_2|\leq |y|\}}(y)+\frac{1}{|y|^{n+\alpha-\gamma}}{\bf 1}_{\{|y|<|x_1-x_2|\}}(y)\bigg)\\[1ex]
&\quad+K\frac{1}{|y|^{n+\alpha}}\| v\|_{1+\gamma'}|x_1-x_2|^{\gamma-\alpha}{\bf 1}_{\{|y| {>} 6\e\}}(y). 
\end{align*}
Therefore, if $\e$ is sufficiently small, then
 \begin{equation}\label{96e}
 \int_{\R^n}C L_2(y)\, dy\leq \Big(\frac{\nu}{2} \|\pi_j^\e v\|_{1+\gamma}+K\| v\|_{1+\gamma'}\Big)|x_1-x_2|^{\gamma-\alpha}
 \end{equation}
 for all~${\tau\in[0,1]},$ $v\in {\rm h}^{1+\gamma}(\bT^n)$, and $0\leq j\leq N$.

Finally, \eqref{f1} (with $\gamma=\gamma'$), the inequality $\beta>\gamma-\alpha$,  and the mean value theorem imply that
\[
 L_3(y)\leq C\|\pi_j^\e v\|_{1+\gamma'} \bigg(\frac{1}{|y|^{n+\alpha}}{\bf 1}_{\{|y|\geq1\}}(y)+\frac{1}{|y|^{n+\alpha-\gamma'}}{\bf 1}_{\{|y|<1\}}(y)\bigg)|x_1-x_2|^{\gamma-\alpha},
\]
hence
\begin{equation}\label{96f}
 \int_{\R^n}  C L_3(y)\, dy\leq  K\| v\|_{1+\gamma'} |x_1-x_2|^{\gamma-\alpha}
 \end{equation}
  for all~${\tau\in[0,1]},$ $v\in {\rm h}^{1+\gamma}(\bT^n)$, and $0\leq j\leq N$.
  
  The estimate \eqref{Dj2} follows now from \eqref{96a}-\eqref{96f}. 
\end{proof}  

It remains to note that  Proposition~\ref{P:1} is established:

\begin{proof}[Proof of Proposition~\ref{P:1}] The assertion follows from Lemma~\ref{L:4P1} and Lemma~\ref{L:5P1}.
\end{proof}

 \subsection{Estimates for the Fourier Multipliers}\label{Sec:3.2}
 In Proposition~\ref{P:1} we have locally approximated the operator $\Phi(\tau u)$ with $\tau\in[0,1]$ and $u\in {\rm h}^{1+\beta}(\bT^n)$ 
 by operators $\delta A^a$ defined in~\eqref{fourmult}  with $\delta\in[1,\eta]$ and  $a\in\R^n$ satisfying $|a|\leq \eta$,  where 
  \begin{equation}\label{etaa}
   \eta:=\big\|(1+|\nabla u|^2)^{1/2}\big\|_0\geq1 .
\end{equation} 
We  first show that these operators are  Fourier multipliers and then  establish in Proposition~\ref{A} a fundamental  estimate for their resolvents.

To start  with,  given $a\in\R^n$ and $k\in\Z^n$, we  use  Fubini's theorem to compute
\begin{align}
&\int_{[-\pi,\pi]^n} \bigg(\int_{\R^n} \frac{\delta_{[x,y]}v-y\cdot\nabla v(x-y)}{|y|^{n+1+\alpha}}\frac{1}{\big[1+\big(|y\cdot a|/|y|)^2\big]^{\tfrac{n+1+\alpha}{2}}}dy\bigg)e^{-ik\cdot x}\, dx\nonumber\\[1ex]
& =\int_{\R^n}\frac{1}{|y|^{n+1+\alpha}}\bigg(\int_{[-\pi,\pi]^n}\big(\delta_{[x,y]}v-y\cdot\nabla v(x-y)\big) e^{-ik\cdot x}\, dx\bigg)\frac{1}{\big[1+\big(|y\cdot a|/|y|)^2\big]^{\tfrac{n+1+\alpha}{2}}}\, dy\nonumber\\[1ex]
&=(2\pi)^{n}\wh v(k)\int_{\R^n}\frac{1-e^{-iy\cdot k}-iy\cdot k e^{-iy\cdot k}}{|y|^{n+1+\alpha}} \frac{1}{\big[1+\big(|y\cdot a|/|y|)^2\big]^{\tfrac{n+1+\alpha}{2}}}\, dy\nonumber\\[1ex]
&=(2\pi)^{n}\wh v(k)\int_{\R^n}\frac{1-\cos(y\cdot k)-y\cdot k\sin (y\cdot k) }{|y|^{n+1+\alpha}}\frac{1}{\big[1+\big(|y\cdot a|/|y|)^2\big]^{\tfrac{n+1+\alpha}{2}}} \, dy,\label{A1}
\end{align}
where we used that
\[
\int_{\R^n}\frac{\sin(y\cdot k)-y\cdot k\cos (y\cdot k) }{|y|^{n+1+\alpha}} \frac{1}{\big[1+\big(|y\cdot a|/|y|)^2\big]^{\tfrac{n+1+\alpha}{2}}}\, dy=0
\]
 since the integrand   is an odd function.
 {The Fourier coefficients $\wh v(k)$, $k\in\Z^n$,  of the periodic function $v$ are defined in Appendix~\ref{Sec:A}.}
Taking into account that $\nabla (1-\cos(y\cdot k))=k\sin(y\cdot k)$, an application of  Gauss' theorem leads us to
\begin{align*}
&\int_{\R^n}\frac{ y\cdot k\sin (y\cdot k) }{|y|^{n+1+\alpha}}\frac{1}{\big[1+\big(|y\cdot a|/|y|)^2\big]^{\tfrac{n+1+\alpha}{2}}} \, dy\nonumber\\[1ex]
&=-\int_{\R^n}(1-\cos(y\cdot k)){\rm div\,}\bigg(\frac{ y  }{|y|^{n+1+\alpha}}\frac{1}{\big[1+\big(|y\cdot a|/|y|)^2\big]^{\tfrac{n+1+\alpha}{2}}}\bigg) \, dy 
\end{align*}
\begin{align}
&=-\int_{\R^n}(1-\cos(y\cdot k)){\rm div\,}\bigg(\frac{ y  }{|y|^{n+1+\alpha}}\Bigg)\frac{1}{\big[1+\big(|y\cdot a|/|y|)^2\big]^{\tfrac{n+1+\alpha}{2}}} \, dy\nonumber\\[1ex]
&=(1+\alpha)\int_{\R^n} \frac{ 1-\cos(y\cdot k)   }{|y|^{n+1+\alpha}} \frac{1}{\big[1+\big(|y\cdot a|/|y|)^2\big]^{\tfrac{n+1+\alpha}{2}}}  \, dy,\label{A2}
\end{align} 
since 
\[
y\cdot\nabla\bigg(\frac{ |y\cdot a|}{|y|}\bigg)^2=0\qquad\text{in $\R^n\setminus\{0\}$.}
\]
Consequently, we derive from \eqref{fourmult}, \eqref{A1}, and \eqref{A2} that
\begin{align}\label{fmult}
\wh{A^a[v]}(k)=m_a(k)\wh v(k),\qquad k\in\Z^n;
\end{align}
that is, $A^a$  defined in \eqref{fourmult} is for $a\in\R^n$  indeed a Fourier multiplier and its symbol is given by
\begin{equation}\label{symbol}
m_a(k):=-2 \int_{\R^n} \frac{ 1-\cos(y\cdot k)   }{|y|^{n+1+\alpha}} \frac{1}{\big[1+\big(|y\cdot a|/|y|)^2\big]^{\tfrac{n+1+\alpha}{2}}}  \, dy,\qquad k\in\Z^n.
\end{equation}

 Our next purpose is to establish the following fundamental resolvent estimate for the Fourier multiplier $A^a$:

{\begin{prop}\label{A}
Given $\eta\ge 1,$ there exists a constant $\kappa=\kappa(\eta)\geq1$ such that
\begin{equation}\label{uniest}
\kappa\|(\lambda-\delta A^a)[v]\|_{\gamma-\alpha}\geq |\lambda|\,\|v\|_{\gamma-\alpha}+\|v\|_{1+\gamma}
\end{equation}
for all $\delta\in[1,\eta]$, $a\in\R^n$ with $|a|\leq \eta,$ 
  $\re\lambda\geq 1$, and $v\in {\rm h}^{1+{\gamma}}(\bT^n)$.
\end{prop}

Proposition~\ref{A} will be a consequence of Lemma~\ref{L:5} and Lemma~\ref{L:6-1}.}
  As a preliminary step, we extend the symbol  $m_a$ defined in \eqref{symbol} to   the whole $\R^n$ by setting 
\begin{equation}\label{symbol2}
m_a(x):=-2 \int_{\R^n} \frac{ 1-\cos(y\cdot x)   }{|y|^{n+1+\alpha}} \frac{1}{\big[1+\big(|y\cdot a|/|y|)^2\big]^{\tfrac{n+1+\alpha}{2}}}  \, dy,\qquad x\in\R^n.
\end{equation}
In  Lemma~\ref{L:4} we  now show that  $m_a$ is smooth on $\R^n\setminus\{0\}$ and in  Lemma~\ref{L:5} we  then establish some estimates  related to $m_a$ which are used in the proof of \eqref{uniest}.
Note that it is not  obvious at first sight that $m_a$  is smooth in~${\R^n\setminus\{0\}}$.
If $a=0$,  this function coincides, up to a multiplicative negative constant, with $|x|^{1+\alpha}$, see the proof of Lemma~\ref{L:4} below. However, if $a\neq 0$, it is not clear 
whether a similar explicit formula  can be derived for  $m_a$.

\begin{lemma}\label{L:4}
  Given $a\in \R^n$, the function $m_a:\R^n\to\R$ defined in \eqref{symbol2} is smooth in~${\R^n\setminus\{0\}}$.
\end{lemma}
\begin{proof} Let $\{e_1,\ldots,e_n\}$  be the canonical basis of $\R^n$.
In order to prove that $m_a$ is smooth in~$\R^n\setminus (e_n[0,\infty))$,
we denote for each given $x\in\R^n\setminus (e_n[0,\infty))$ by  $\mathsf{H}_n=\mathsf{H}_n(x) \in \R^{n\times n}$ the Householder transformation
\[
\mathsf{H}_n:=I_n-2w w^\top,
\]
where 
$$
w:=w(x):=\frac{x/|x|-e_n}{|x/|x|-e_n|}\in\R^n
$$ 
 and $I_n\in \R^{n\times n}$ is the identity matrix.
 Then $\mathsf{H}_n$ is symmetric, orthogonal, and $\mathsf{H}_nx=|x|e_n$.
 Let further 
 \[
 {f_n:V_n:=(0,\infty)\times(0,2\pi)\times(0,\pi)^{n-2}\to  \R^n},
 \] 
 with $f_n=f_n(r,\varphi,\vartheta)=rg_n(\varphi,\theta)$, $\vartheta=(\vartheta_1,\ldots,\vartheta_{n-2})$, and $|g_n(\varphi,\theta)|=1$
 be the standard $n$-dimensional  polar coordinates transformation, see e.g~\cite{AE3}.
   Write $g_n=(g_{n}^{1},\ldots, g_{n}^{n})$.
  Then, changing  variables according to $y=\mathsf{H}_nf_n(r,\varphi,\vartheta)$ and using the orthogonality and symmetry of $\mathsf{H}_n$, we have
 \begin{align*}
 &\int_{\R^n} \frac{ 1-\cos(y\cdot x)   }{|y|^{n+1+\alpha}}  \frac{1}{\big[1+\big(|y\cdot a|/|y|)^2\big]^{\tfrac{n+1+\alpha}{2}}} \, dy   \\[1ex]
 &=\int_{V_n} \frac{1-\cos( f_n(r,\varphi,\vartheta)\cdot \mathsf{H}_nx)}{|  f_n(r,\varphi,\vartheta) |^{n+1+\alpha}} 
  \frac{|\det(\p f_n(r,\varphi,\vartheta))|}{\big[1+\big(|   f_n(r,\varphi,\vartheta\cdot \mathsf{H}_n a)|/|  f_n(r,\varphi,\vartheta)|\big)^2\big]^{\tfrac{n+1+\alpha}{2}}}  \, d(r,\varphi,\vartheta)\\[1ex]
    &=\int_{V_n} \frac{1-\cos( |x|r g_n^n(\varphi,\vartheta)))}{r^{2+\alpha}}
    \frac{\sin(\vartheta_1)\sin^2(\vartheta_2)\ldots \sin^{n-2}(\vartheta_{n-2})}{\big[1+|    g_n(\varphi,\vartheta)\cdot \mathsf{H}_n a|^2\big]^{\tfrac{n+1+\alpha}{2}}}\, d(r,\varphi,\vartheta)\\[1ex]
     &=|x|^{1+\alpha}\int_{V_n} \frac{1-\cos( \tau g_n^n(\varphi,\vartheta))}{\tau^{2+\alpha}}
     \frac{\sin(\vartheta_1)\sin^2(\vartheta_2)\ldots \sin^{n-2}(\vartheta_{n-2})}{\big[1+|   g_n(\varphi,\vartheta)\cdot \mathsf{H}_n a|^2\big]^{\tfrac{n+1+\alpha}{2}}}\, d(\tau,\varphi,\vartheta),
 \end{align*}
 hence 
 \[
m_a(x)=|x|^{1+\alpha}p_n(x),\qquad x\in \R^n\setminus (e_n[0,\infty)),
 \]
 where $p_n:\R^n\setminus (e_n[0,\infty))\to (-\infty,0)$  is defined by the relation
 \begin{equation}\label{repr}
p_n(x):=-2 \int_{V_n} \frac{1-\cos( \tau g_n^n(\varphi,\vartheta))}{\tau^{2+\alpha}}
    \frac{\sin(\vartheta_1)\sin^2(\vartheta_2)\ldots \sin^{n-2}(\vartheta_{n-2})}{\big[1+|\mathsf{H}_n(x) a\cdot   g_n(\varphi,\vartheta)|^2\big]^{\tfrac{n+1+\alpha}{2}}}\, d(\tau,\varphi,\vartheta) 
 \end{equation}
 with $\mathsf{H}_n=\mathsf{H}_n(x)=(h_{ij}(x))_{1\leq i,\,j\leq n}$ and
 \begin{equation}\label{hij}
h_{ij}(x):=\delta_{ij}-\Big(\frac{x_i}{|x|}-\delta_{in}\Big) \Big(\frac{x_j}{|x|}-\delta_{jn}\Big)\Big(1-\frac{x_n}{|x|}\Big)^{-1}.
 \end{equation}
Clearly, all  functions $h_{ij}$, $1\leq i,\,j\leq n,$ are smooth in $\R^n\setminus (e_n[0,\infty))$ and their partial derivatives are all bounded on compact subsets of $\R^n\setminus (e_n[0,\infty))$.
 A direct application of the theorem on differentiation under the integral sign allows one to conclude that~$p_n$, and hence also $m_a$, is smooth in $\R^n\setminus (e_n[0,\infty))$.

Arguing similarly  as above, but using  instead of $\mathsf{H}_n$ the Householder transformation ${\mathsf{H}_1=\mathsf{H}_1(x)}$  which maps~${x\in\R^n\setminus (e_1[0,\infty))}$
to~$|x|e_1$, we may deduce that $m_a$ is smooth also in~${\R^n\setminus (e_1[0,\infty))}$ and the desired claim follows.
\end{proof}

In view of Lemma~\ref{L:4}, the function $P_a:\R^n\setminus\{0\}\to\R$, defined for fixed $a\in\R^n$ by
\begin{equation}\label{defp}
  P_a (x):=\frac{m_a(x)}{|x|^{1 +\alpha}},\qquad x\in\R^n\setminus\{0\},
\end{equation}
 is smooth, and we infer from the definition \eqref{symbol2} that $P_a(x)=P_a(-x)$ for all $x\in\R^n\setminus\{0\}$.
 Lemma~\ref{L:5} below provides some useful estimates on $P_a$. 
  Before stating this result, we note that if $\cO\subset \R^n$ is an open set and $f:\cO\to\R$ is a positive function,
 then for each $ {0\neq\mu}\in\N^n$, we have
  \begin{equation}\label{faadibruno}
\p^\mu   \Big(\frac{1}{f}\Big)=\sum_{k=1}^{|\mu|}\frac{1}{f^{k+1}}
\underset{\substack{
\{\beta_1,\ldots,\beta_k\}\subset\N^n\\
\beta_1+\ldots+\beta_k=\mu
}}\sum
c_{\beta_1\ldots\beta_k}\p^{\beta_1}f\cdot\ldots \cdot\p^{\beta_k}f,
 \end{equation}
 with constants $c_{\beta_1\ldots\beta_k}\in\Z$ independent of $f$.
 The formula \eqref{faadibruno} is obtained by using a standard induction argument.
  {In Lemma~\ref{L:5}-Lemma~\ref{L:6} below we use the notation  
 \begin{equation}\label{NN}
 N:=[n/2]+1, 
 \end{equation}
 where $[n/2]$ is the integer part of $n/2$.}
   
\begin{lemma}\label{L:5}
Let $P_a$ be defined in \eqref{defp}.
Given $\eta\geq1$, there exists a constant \mbox{$M=M(\eta)\geq 1$} such that for all~${a\in\R^n}$ with $|a|\leq \eta$
we have
\begin{equation}\label{estfour}
M^{-1}\leq  \inf_{x\in\R^n\setminus\{0\}}|P_a(x)|\qquad\text{and}\qquad \sup_{x\in\R^n\setminus\{0\}}  |x|^{|\mu|}|\p^\mu P_a(x)|\leq M,
\end{equation}
whenever $\mu\in\N^n$  satisfies $|\mu|\leq N$.
\end{lemma}
 \begin{proof}
  Since $P_a$ is smooth and $P_a(x)=P_a(-x)$ for all $x\in\R^n\setminus\{0\}$, it suffices to estimate~$P_a$ and its  partial derivatives   in the set~$\{x_n<0\}$, in which case we also have
  $P_a=p_n$, see~\eqref{repr}.
  It readily follows from \eqref{repr} that, if  $x\in \R^n\setminus (e_n[0,\infty))$, then
\begin{equation*} 
|p_n(x)|\geq \frac{2}{\big[1+\eta^{2}\big]^{\tfrac{n+1+\alpha}{2}}} \int_{V_n} \frac{1-\cos( \tau g_n^n(\varphi,\vartheta))}{\tau^{2+\alpha}}
     \sin(\vartheta_1)\sin^2(\vartheta_2)\ldots \sin^{n-2}(\vartheta_{n-2}) \, d(\tau,\varphi,\vartheta).
 \end{equation*} 
 This proves the first estimate in \eqref{estfour}.
 
  Observing for $x\in\{x_n<0\}$ that
  \[
  1-\frac{x_n}{|x|}\geq 1,
  \]
 it follows that the  functions $h_{ij}$, $1\leq i,\,j\leq n,$ defined in \eqref{hij} are smooth in $\{x_n<0\}$ and 
\[
\max_{|\mu|\leq N} \sup_{\{x_n<0\}}|x|^{|\mu|}|\p^\mu h_{ij}(x)|<\infty \qquad\text{for all $1\leq i,\, j\leq n$}. 
 \]
  A direct application of the theorem on the 
  differentiation under the integral sign and of~\eqref{faadibruno} enables us to deduce also the  second estimate in \eqref{estfour}.
 \end{proof}

 With these preparatory results we may now establish  the resolvent estimate on $\delta A^a$ stated in~Proposition~\ref{A}.  
 For this purpose, let $\eta\ge 1$ be as in~\eqref{etaa} and consider  $\delta\in[1,\eta]$,  $a\in\R^n$ with~$|a|\leq \eta$, and  $\re\lambda\geq1$. We introduce $R(\lambda)$ as the Fourier multiplier
 \begin{equation}\label{RL}
\wh{R(\lambda)[v]}(k):=\frac{ \mathfrak{m}_{1+\alpha}(k)}{\lambda-\delta m_a(k)} \wh v (k),\qquad k\in\Z^n,
 \end{equation}
 with   $m_a$ being the Fourier symbol of $A^a$ defined in \eqref{symbol2} and $\mathfrak{m}_{1+\alpha}$ 
is a smooth function defined  on  $\R^n$ such that $\mathfrak{m}_{1+\alpha}(x)=1$ for $|x|\leq1/2$ and $\mathfrak{m}_{1+\alpha}(x)=|x|^{1+\alpha}$ for $|x|\geq 1$. 
 Using the Mikhlin-H\"ormander multiplier theorem  we shall show below that $R(\lambda)\in\kL({\rm h}^{\gamma-\alpha}(\bT^n))$ and 
 that there exists $\kappa=\kappa(\eta)$ such that 
\begin{equation}\label{bound}
\|R(\lambda)\|_{\kL({\rm h}^{\gamma-\alpha}(\bT^n))}\leq \kappa ,
\end{equation}
whenever  $\delta\in[1,\eta]$,  $a\in\R^n$ with~$|a|\leq \eta$, and  $\re\lambda\geq1$. We also note that the operator~$I_{-1-\alpha}$ from Lemma~\ref{L:A1} with Fourier symbol $1/\mathfrak{m}_{1+\alpha}$ satisfies
\begin{equation}\label{xy}
I_{-1-\alpha}\in \kL({\rm h}^{\gamma-\alpha}(\bT^n),{\rm h}^{1+\gamma}(\bT^n)).
\end{equation}
Indeed, \eqref{xy} readily follows from Lemma~\ref{L:A1},  the fact that one may identify ${\rm C}^s(\bT^n)$
with the Besov space   $B^{s}_{\infty,\infty}(\bT^n)$  for~${s>0}$ with $s\not\in\N$, e.g. see \cite{ST87}, and the fact that $I_{-1-\alpha}$ maps~${{\rm C}^\infty(\bT^n)}$ to ${\rm C}^\infty(\bT^n)$. Therefore, we infer from \eqref{bound} and \eqref{xy} that
\begin{equation}\label{x}
(\lambda-\delta A^a)^{-1}=I_{-1-\alpha}R(\lambda)\in \kL({\rm h}^{\gamma-\alpha}(\bT^n), {\rm h}^{1+\gamma}(\bT^n)).
\end{equation}  
This is the content of the following lemma.

\begin{lemma}\label{L:6-1} 
 Given $\eta\geq1$, there exists a constant $\kappa=\kappa(\eta)\geq1$  such that if~${\delta\in[1,\eta]}$,  $a\in\R^n$ with~$|a|\leq \eta$,  $\re\lambda\geq1$, 
 and $v\in{\rm h}^{1+\gamma}(\bT^n)$, then
\begin{equation}\label{uniestt0}
\kappa\|(\lambda-\delta A^a)[v]\|_{\gamma-\alpha}\geq \|v\|_{1+ \gamma}.
\end{equation}
Moreover, $\lambda-\delta A^a\in\kL({\rm h}^{1+ \gamma}(\bT^n), {\rm h}^{\gamma-\alpha}(\bT^n))$ is an isomorphism.
 \end{lemma}
\begin{proof}
As pointed out above, \eqref{uniestt0} follows from~\eqref{bound}, \eqref{xy}, and~\eqref{x}. It thus remains to establish \eqref{bound}, 
for which the  Mikhlin-H\"ormander multiplier theorem \cite[Section~3.6.3/Remark~3]{ST87} entails that it suffices to show that
 \begin{equation}\label{FourE2}
 \sup_{|\beta|\leq N}\sup_{x\in\R^n\setminus\{0\}}|x|^{|\beta|}\Big|\p^\beta \Big(\frac{|x|^{1+\alpha}}{\lambda-\delta m_a(x)}\Big)\Big|\leq \kappa,
 \end{equation}
uniformly with respect to~${\delta\in[1,\eta]}$,  $a\in\R^n$ with~$|a|\leq \eta$, and $\re\lambda\geq1$,  {with~$N$ given by~\eqref{NN}.}   
 Let thus  {$\beta\in\N^n$} satisfy $|\beta|\leq N$. 
 In view of Leibniz' rule, it remains to estimate terms of the form
\[
|x|^\beta \big|\p^{\beta-\mu}|x|^{1+\alpha}\big|\Big|\p^\mu \Big(\frac{1}{\lambda-\delta m_a(x)}\Big)\Big|
\] 
with $\mu\in\N^n$ such that  $\mu\leq \beta$.
 To this end, first note that 
\begin{equation}\label{F11}
 \big|\p^{\beta-\mu}|x|^{1+\alpha}\big|\leq C|x|^{1+\alpha+|\mu|- {|\beta|}},
\end{equation}
where $C=C(N)$. 
Moreover, it follows from \eqref{faadibruno} that
 \begin{equation}\label{faa2}
\p^\mu \big((\lambda-\delta m_a)^{-1}\big)=\sum_{k=1}^{|\mu|}\frac{\delta^k }{(\lambda-\delta m_a)^{k+1}}
\underset{\substack{
\{\beta_1,\ldots,\beta_k\}\subset\N^n\\
\beta_1+\ldots+\beta_k=\mu
}}\sum
c_{\beta_1\ldots\beta_k}\p^{\beta_1}m_a\cdot\ldots \cdot\p^{\beta_k}m_a,
 \end{equation}
 with constants $c_{\beta_1\ldots\beta_k}\in\Z$ independent of $\lambda$, $\delta$, and $m_a$.
We note that $m_a(x)<0$ for~$x\neq0$.
   Therefore, given $\delta\in[1,\eta]$, $\re\lambda\geq 1$, and $x\in\R^n\setminus\{0\}$, we infer from \eqref{defp} and  Lemma~\ref{L:5} that
 \begin{equation}\label{bddb}
 \begin{aligned}
 |\lambda-\delta m_a(x)|&=\sqrt{(\re\lambda-\delta m_a(x))^2+(\im \lambda)^2}\geq\max\{|\lambda|, |m_a(x)|\}\\[1ex]
 &\geq \max\{|\lambda|, M^{-1}|x|^{1+\alpha}\},
 \end{aligned}
 \end{equation}
 with $M=M(\eta)\geq 1$, and therefore
 \begin{equation}\label{bdd1}
 \frac{1}{|\lambda-\delta m_a(x)|^{k+1}}\leq \frac{M^{k+1}}{|x|^{(k+1)(1+\alpha)}}, \qquad  k\ge 0.
 \end{equation}
 Leibniz' rule, \eqref{F11},  and Lemma~\ref{L:5} ensure the existence of a constant~${C=C(\eta)}$ such that 
 for all $x\in\R^n\setminus\{0\}$ we have
\begin{equation}\label{bddb2}
|(\p^{\beta_1}m_a\cdot\ldots \cdot\p^{\beta_k}m_a)(x)|\leq C|x|^{1+\alpha-|\beta_1| }\ldots |x|^{1+\alpha-|\beta_k| }=C|x|^{k(1+\alpha)-|\mu|}.
\end{equation}
The estimate \eqref{FourE2} is now  a straightforward consequence of \eqref{F11}, \eqref{bdd1}, and \eqref{bddb2}.
\end{proof}

To complete the proof of Proposition~\ref{A}, we establish the following resolvent estimate:
 \begin{lemma}\label{L:6} 
 Given $\eta\geq1$, there exists a constant $\kappa=\kappa(\eta)\geq1$  such that 
\begin{equation}\label{uniestt1}
\|(\lambda-\delta A^a)^{-1}\|_{\kL({\rm h}^{\gamma-\alpha}(\bT^n))}\leq \frac{\kappa}{|\lambda|}
\end{equation}
 for all~${\delta\in[1,\eta]}$,  $a\in\R^n$ with~$|a|\leq \eta$, and $\re\lambda\geq1$.
 \end{lemma}

\begin{proof}
 Since $(\lambda-\delta A^a)^{-1}=I_{-1-\alpha}R(\lambda)$,  the Mikhlin-H\"ormander multiplier theorem \cite[Section~3.6.3/Remark~3]{ST87} implies that we only have to show that
 \begin{equation}\label{FourE1}
 \sup_{|\mu|\leq N}\sup_{x\in\R^n\setminus\{0\}}|x|^{|\mu|}\Big|\p^\mu \big((\lambda-\delta m_a)^{-1}\big)(x)\Big|\leq\frac{\kappa}{|\lambda|}
 \end{equation}
  uniformly in~${\delta\in[1,\eta]}$,  $a\in\R^n$ with~$|a|\leq \eta$, and~${\re\lambda\geq1}$.  
 
Invoking  \eqref{bddb}, there exists a constant $M=M(\eta)\geq 1$ such that  for  all ${\delta\in[1,\eta]}$,~${\re\lambda\geq 1}$, and~${x\in\R^n\setminus\{0\}}$, we have
 \[
 \frac{1}{|\lambda-\delta m_a(x)|^{k+1}}\leq \frac{ 1}{|\lambda|}\frac{M^k}{|x|^{k(1+\alpha)}}, \qquad  k\ge 0.
 \]
The  desired claim \eqref{FourE1} follows now from the latter estimate,  \eqref{faa2}, and \eqref{bddb2}.
 \end{proof}

\begin{proof}[Proof of Proposition~\ref{A}] This is now a consequence of  Lemma~\ref{L:5} and Lemma~\ref{L:6}.
\end{proof}

\subsection{ Proof of Theorem~\ref{T:GP}}\label{Sec:3.3}
Having established Proposition~\ref{P:1} and~Proposition~\ref{A}, we turn to the proof of Theorem~\ref{T:GP}.
\begin{proof}[Proof of Theorem~\ref{T:GP}]
Let  {$\gamma'\in(\max\{\alpha,\beta\},\gamma)$}, let $\eta\geq1$ be as defined in \eqref{etaa}, and  let~${\kappa\geq1}$ be chosen as in Proposition~\ref{A}.
 Defining $\nu:=(2\kappa)^{-1}$, Proposition~\ref{P:1} ensures there exist~${\e\in(0,1)}$ and  a positive constant~${K=K(\e)}$ with the property that for all~${\tau\in[0,1]},$ $v\in {\rm h}^{1+\gamma}(\bT^n)$, 
 and $0\leq j\leq N$ we have
 \begin{equation}\label{Dj21p}
  2\kappa\big\|\pi_j^\e \Phi(\tau u)[v]- (1+\tau^2|\nabla u|^2)^{1/2}(x_j^\e)A^{\tau \nabla u(x_j^\e)}[\pi_j^\e v]\big\|_{\gamma-\alpha}\leq   \|\pi_j^\e v\|_{1+\gamma}+2\kappa K\|  v\|_{1+\gamma'}.
 \end{equation}
 {We now} infer from Proposition~\ref{A}  that
 for all~${\tau\in[0,1]},$ $v\in {\rm h}^{1+\gamma}(\bT^n)$, $0\leq j\leq N$, and $\re\lambda\geq 1$   
\begin{equation}\label{uniestp}
2\kappa\|(\lambda-(1+\tau^2|\nabla u|^2)^{1/2}(x_j^\e)A^{\tau \nabla u(x_j^\e)})[\pi_j^\e v]\|_{\gamma-\alpha}\geq 2|\lambda|\,\|\pi_j^\e v\|_{\gamma-\alpha}+2\|\pi_j^\e v\|_{1+\gamma}.
\end{equation}
Combining \eqref{Dj21p} and \eqref{uniestp}, we get
 \begin{align*}
   2\kappa\|\pi_j^\e(\lambda-\Phi(\tau u))[v]\|_{ \gamma-\alpha}\geq& 2\kappa\|(\lambda-(1+\tau^2|\nabla u|^2)^{1/2}(x_j^\e)A^{\tau \nabla u(x_j^\e)})[\pi_j^\e v]\|_{\gamma-\alpha}\\[1ex]
   &\quad -2\kappa\big\|\pi_j^\e \Phi(\tau u)[v]- (1+\tau^2|\nabla u|^2)^{1/2}(x_j^\e)A^{\tau \nabla u(x_j^\e)}[\pi_j^\e v]\big\|_{\gamma-\alpha}\\[1ex]
   \geq& 2|\lambda|\,\|\pi^\e_j v\|_{\gamma-\alpha}+ \|\pi^\e_j v\|_{1+\gamma}-2\kappa K\|  v\|_{1+\gamma'}.
 \end{align*}
 Summing  over $j$, we  deduce from \eqref{eqnorm}, Young's inequality, and the interpolation property~\eqref{interpol}
that there exist constants  $\kappa_1\geq1$  and~$\omega_1\geq1 $ such that 
  \begin{align}\label{KDED}
   \kappa_1\|(\lambda-\Phi(\tau u))[v]\|_{\gamma-\alpha}\geq  |\lambda| \|v\|_{\gamma-\alpha}+ \| v\|_{1+\gamma}
 \end{align}
for all~${\tau\in[0,1]},$ $v\in {\rm h}^{1+\gamma}(\bT^n)$,  and $\re\lambda\geq \omega_1.$ 

Moreover,  Lemma~\ref{L:6-1} ensures that   $$\omega_1-\Phi(0)=\omega_1-A^0 \in \kL({\rm h}^{1+\gamma}(\bT^n), {\rm h}^{\gamma-\alpha}(\bT^n))$$ is an isomorphism.
This property, the method of continuity \cite[Proposition I.1.1.1]{Am95}, and~\eqref{KDED}   imply now that also
$\omega_1-\Phi(u)\in\kL({\rm h}^{1+\gamma}(\bT^n), {\rm h}^{\gamma-\alpha}(\bT^n))$ is an isomorphism.
 Consequently, together with  \eqref{KDED} (with $\tau=1$) we infer from \cite[Chapter~I]{Am95} that indeed
\[
-\Phi(u)\in\mathcal{H}({\rm h}^{1+\gamma}(\bT^n), {\rm h}^{ \gamma-\alpha}(\bT^n))
\]
as claimed.
\end{proof}

 \subsection{ Proof of Theorem~\ref{MT1}}\label{Sec:3.4}
 
 We are now in a position to establish our first  main  result stated in Theorem~\ref{MT1}.
 \begin{proof}[Proof of Theorem~\ref{MT1}]
   {Let  $\beta'\in(0,\beta)$ be chosen such that $\gamma<\alpha+\beta'$.} Combining  Proposition~\ref{P1} and Theorem~\ref{T:GP}
 we have
 \begin{equation}\label{property2}
-\Phi\in{\rm C}^\infty\big({\rm h}^{1+\beta'}(\bT^n),\mathcal{H}({\rm h}^{1+\gamma}(\bT^n), {\rm h}^{ \gamma-\alpha}(\bT^n))\big).
\end{equation}
Moreover, setting 
\[
\theta':=1-\frac{\gamma-\beta'}{1+\alpha}\qquad\text{and}\qquad \theta:=1-\frac{\gamma-\beta}{1+\alpha},
\]
we infer from \eqref{interpol} that
\begin{equation*}
({\rm h}^{\gamma-\alpha}(\bT^n), {\rm h}^{1+\gamma}(\bT^n))_{\theta',\infty}^0={\rm h}^ {1+\beta'}(\bT^n)\qquad\text{and}\qquad
({\rm h}^{\gamma-\alpha}(\bT^n), {\rm h}^{1+\gamma}(\bT^n))_{\theta,\infty}^0={\rm h}^ {1+\beta}(\bT^n).
\end{equation*}
The property~\eqref{property2} and the latter  interpolation relations ensure that we may  apply \cite[Theorem 1.1]{MW20}  to the evolution problem~\eqref{QAP}.
Hence, given~${u_0\in {\rm h}^ {1+\beta}(\bT^n)}$, there exists a unique maximal classical solution  $u= u(\,\cdot\, ; u_0)$ to~\eqref{QAP} such that
\begin{equation}\label{reggg}
 u\in {\rm C}([0,T^+),{\rm h}^{1+\beta}(\bT^n))\cap {\rm C}((0,T^+),{\rm h}^{1+\gamma}(\bT^n))\cap {\rm C}^1((0,T^+),{\rm h}^{\gamma-\alpha}(\bT^n))
  \end{equation}
   and  
  \begin{equation*} 
u\in   {\rm C}^{\theta-\theta'}([0,T^+), {\rm h}^{1+\beta'}(\bT^n)),
  \end{equation*}
  where $T^+=T^+(u_0)\in(0,\infty]$ denotes the maximal existence time.
  Moreover, the mapping~${[(t,u_0)\mapsto u(t;u_0)]}$ defines a semiflow on ${\rm h}^{1+\beta}(\bT^n)$ which is smooth in the open set
  \[
\{(t,u_0)\,:\, u_0\in {\rm h}^{1+\beta}(\bT^n),\, 0<t<T^+(u_0)\}\subset \R\times {\rm h}^{1+\beta}(\bT^n).  
  \]
  As shown in Lemma~\ref{L:US} below, the solution is  actually unique in the larger set of functions satisfying only \eqref{reggg}.

Since the embedding ${\rm h}^{1+\gamma}(\bT^n)\hookrightarrow {\rm h}^{ \gamma-\alpha}(\bT^n)$ is compact,  Theorem~\ref{MT1}~(i) follows directly from \cite[Theorem~1.1~(iv)$(\beta)$]{MW20}.

We next establish the parabolic smoothing property~\eqref{eq:fg}. 
This may be shown by using a parameter trick  employed also  in other settings, see \cite{An90, ES96, PSS15}. 
The arguments are more or less identical to those presented in \cite[Theorem 1.3]{MBV19} and we include them here merely for the reader's ease.

To start  with, we fix a maximal solution $u=u(\cdot; u_0)$ to \eqref{QAP}. 
Given   $\lambda_1\in (0,\infty)$ and~${\lambda_2\in\R^n}$, we set $\lambda:=(\lambda_1,\lambda_2)$ and define  
\[
u_{\lambda}(t,x):=u(\lambda_1 t,x+\lambda_2 t), \qquad   x\in\R^n, \quad 0\leq t<T^+(u_0)/\lambda_1.
\]
Using the density of ${\rm C}^\infty(\bT^n) $ in the little H\"older spaces, it is possible to show that 
\begin{equation}\label{reggggg}
 u_\lambda\in {\rm C}([0,T^+_\lambda),{\rm h}^{1+\beta}(\bT^n))\cap {\rm C}((0,T^+_\lambda),{\rm h}^{1+\gamma}(\bT^n))\cap {\rm C}^1((0,T^+_\lambda,{\rm h}^{\gamma-\alpha}(\bT^n)), 
  \end{equation}
where $T^+_\lambda:=T^+(u_0)/\lambda_1$, and
 \begin{equation}\label{eqqu}
\frac{du_\lambda}{dt}(t)=\lambda_1\Phi(u_\lambda(t))[u_\lambda(t)]+\lambda_2\cdot\nabla u_\lambda(t),\qquad t\in(0,T^+_\lambda).  
 \end{equation}
Let now $U :=(U_1,U_2,U_3):[0,T^+_\lambda)\to\R^{n+1}\times {\rm h}^{1+\beta}(\bT^n)$ be the function defined by
\[(U_1,U_2,U_3)(t):=(\lambda_1,\lambda_2, u_\lambda(t)).\] 
Then, in view of \eqref{eqqu}, the function $U$ is a solution to the quasilinear evolution problem
\begin{align}\label{QC}
\frac{d U}{dt}(t)= \Psi(U(t))[U(t)],\quad t>0,\qquad U(0)=U_0:=(\lambda,u_0),
\end{align}
where 
$$ \Psi:(0,\infty)\times\R^n\times {\rm h}^{1+\beta}(\bT^n)\to\kL\big(\R^{n+1}\times {\rm h}^{1+\gamma}(\bT^n), \R^{n+1}\times{\rm h}^{ \gamma-\alpha}(\bT^n)\big)
$$ 
is the operator defined by
\begin{align}\label{QO}
  \Psi(U_1,U_2,U_3)[(V_1,V_2,V_3)]:=\big(0,0,U_1\Phi(U_3)[V_3]+U_2\cdot\nabla V_3\big).
\end{align}
Since $\nabla $ is a first order operator and $-\Phi(U_3)\in\mathcal{H}({\rm h}^{1+\gamma}(\bT^n),{\rm h}^{ \gamma-\alpha}(\bT^n))$, we infer from \eqref{property}  and \cite[Corollary I.1.6.3]{Am95} that 
\begin{align*}
 -\Psi\in {\rm C}^\infty\big((0,\infty)\times\R^n\times {\rm h}^{1+\beta}(\bT^n),\mathcal{H}(\R^{n+1}\times {\rm h}^{1+\gamma}(\bT^n), \R^{n+1}\times{\rm h}^{ \gamma-\alpha}(\bT^n))\big). 
\end{align*}

Hence, we may argue as above to conclude that  the quasilinear parabolic problem \eqref{QC} has for each 
 $U_0=(\lambda,u_0)\in (0,\infty)\times\R^n\times {\rm h}^{1+\beta}(\bT^n)$  a unique maximal solution~${U=U(\cdot;U_0)}$ with $(U_1,U_2)=\lambda_2$ and
\[ U_3\in {\rm C}([0,T^+(U_0)),{\rm h}^{1+\beta}(\bT^n))\cap {\rm C}((0,T^+(U_0)),{\rm h}^{1+\gamma}(\bT^n))\cap {\rm C}^1((0,T^+(U_0)),{\rm h}^{\gamma-\alpha}(\bT^n)) .\]
In addition,
\[
\0:=\{(\lambda,u_0,t)\,:\, t\in(0,T^+(\lambda,u_0))\}
\]
is an open subset of $\R^{n+1}\times {\rm h}^{1+\beta}(\bT^n)\times(0,\infty)$ and 
\[
[(\lambda,u_0,t)\mapsto U(t;(\lambda,u_0))]:\0\to \R^{n+1}\times {\rm h}^{1+\beta}(\bT^n)
\]
is a smooth mapping. Therefore, given  $u_0\in {\rm h}^{1+\beta}(\bT^n),$  we have
\[
\frac{T^+(u_0)}{\lambda_1}=T^+(\lambda,u_0) \qquad\text{for  $\lambda\in(0,\infty)\times\R^{n}.$}
\]

In order to prove that  $u=u(\cdot;u_0)$ is smooth on $(0,T^+(u_0))\times\R^n$, we choose  an arbitrary point $(t_0,x_0)\in (0,T^+(u_0))\times\R^n$.
Let $\bB_\e\subset\R^{n+1}$  be the ball centered at $(1,\ldots,1)$ of radius $\e>0$, where we choose $\e$ sufficiently small to guarantee that
for all $\lambda\in\bB_\e$ we have~${t_0<T^+(\lambda,u_0).}$
In particular, we have $\bB_\e\times\{u_0\}\times\{t_0\}\subset \0 $ and the restriction
\begin{align*}
[\lambda\mapsto u_\lambda(t_0)]:\bB_\e\to {\rm h}^{1+\beta}(\bT^n)
\end{align*}
is  smooth too, hence also  the mapping
\begin{align}\label{11}
[\lambda\mapsto u_\lambda(t_0,x_0)=u(\lambda_1t_0,x_0-t_0+\lambda_2t_0)]:\bB_\e\to \R
\end{align}
is smooth.
Letting $\bB_\delta(t_0,x_0)$ be the ball centered at $(t_0,x_0)$ with sufficiently small radius $\delta$, we note that the function $\varphi:\bB_\delta(t_0,x_0)\to\bB_\e$, defined by
  \begin{align}\label{12}
  \varphi(t,x):= \Big(\frac{t}{t_0},\frac{x-x_0+t_0}{t_0}\Big),
  \end{align}
is well-defined and smooth.
 Composing   the functions defined in \eqref{11} and \eqref{12}, we conclude that 
$$  \big[(t,x)\mapsto  u(t,x)\big]:\bB_\delta (t_0,x_0)\to\R$$
is indeed smooth.  This yields Theorem~\ref{MT1}~(ii).

The proof  for the a priori estimates of Theorem~\ref{MT1}~(iii) is provided in Lemma~\ref{L:9} below.
 \end{proof}

\section{Stability Analysis}\label {Sec:5}
This section is devoted to the proof of Theorem~\ref{MT2}.
Let us emphasize that since each constant is a stationary solution to \eqref{QAP}, the linearization of the right-hand side of \eqref{QAP} 
at any constant solution has zero as an eigenvalue, and this impedes us to 
 applying directly (some version of the) the principle linearized stability  in the context of~\eqref{QAP}.  
In order to establish our stability result presented in Theorem~\ref{MT2}, 
we introduce a volume preserving unknown solving a  quasilinear evolution problem, see~\eqref{NEVS} and~\eqref{QAP'S} below,  to which  
 the quasilinear principle of linearized stability established in \cite[Theorem~1.3]{MW20}  applies.
  {Doing this, we  derive in \eqref{stestS} and \eqref{hgh} estimates  for the volume preserving unknown 
 which are key to prove  the exponential decay estimate \eqref{Exst}.
 
  {We begin this section by showing that the stationary solutions to \eqref{QAP}  coincide with the constant functions, see Lemma~\ref{L:8} below.
In Section~\ref{Ssec:52} we provide the  improved uniqueness statement still missing in the proof of Theorem~\ref{MT1} and some maximum principles for the nonlocal mean curvature flow.
We conclude the section with the proof of the stability result stated in Theorem~\ref{MT2}.

Throughout this section the H\"older exponents $\alpha,\,\beta,\, \gamma$ are assumed to satisfy \eqref{constants}.}

\subsection{Stationary Solutions to \eqref{QAP}}\label{Ssec:51}
Recalling~\eqref{NMC}, it is straightforward to see that each constant is a stationary solution to \eqref{QAP}. 
In order to establish  that these are the only stationary solutions, we use an alternative expression  for the nonlocal mean curvature $H_\alpha(u)$ provided, e.g., in \cite{Fall18, AFW22xx}.
More precisely, using Gauss' formula,  one may reformulate the nonlocal mean curvature $H_\alpha(u)$ for given~ {$u\in {\rm h}^{1+\gamma}(\bT^n)$}  as
\begin{equation}\label{NMC2} 
H_\alpha(u)(x)=-\PV\int_{\R^n}  \frac{1}{|y|^{n+\alpha}}\bigg[ F \Big(  \frac{\delta_{[x,y]}u}{|y|} \Big)-F \Big( - \frac{\delta_{[x,y]}u}{|y|} \Big)\bigg]\, dy ,\qquad x\in\R^n,
 \end{equation}
 where the function $F:\R\to\R$ is  defined by
 \[
F(\xi):=\int_\xi^\infty \frac{d\tau}{(1+\tau^2)^{\tfrac{n+1+\alpha}{2}}},\qquad \xi\in\R.
\]

The following result  identifies the stationary solutions as the constant constants  and  relies on \cite[Lemma 3.6]{Fall18}. 

 \begin{lemma}\label{L:8} Let $u\in{\rm C}^{1+\gamma}(\bT^n)$  be  such that  $H_\alpha(u)=0$.
 Then $u$ is a constant. 
 \end{lemma}

 \begin{proof}
 Let $m:=\min_{\R^n}u$ and assume that $v=u-m$ is not identically zero.
 Then, there exists $x_0\in\R^n$ such that $v(x_0)=\max_{\R^n} v>0$. In view of $H_\alpha(u)=0$ and~\eqref{NMC2} we get
 \begin{equation}\label{NMC2'} 
\PV\int_{\R^n} \frac{1}{|y|^{n+\alpha}}\bigg[ F \Big(  \frac{\delta_{[x_0,y]}v}{|y|} \Big)-F \Big( - \frac{\delta_{[x_0,y]}v}{|y|} \Big)\bigg]\, dy=0.
 \end{equation}
Moreover,  the mean value theorem and  $F'<0$ imply for $y\neq 0$ that
\[
F \Big(  \frac{\delta_{[x_0,y]}v}{|y|} \Big)-F \Big( - \frac{\delta_{[x_0,y]}v}{|y|} \Big)=2\frac{v(x_0)-v(x_0-y)}{|y|}\int_0^1 F'\Big(  s\frac{\delta_{[x_0,y]}v}{|y|}-(1-s) \frac{\delta_{[x_0,y]}v}{|y|} \Big)\, ds\leq0,
\]
the latter expression being (strictly) negative when $x_0-y$ is close to a point where  $v$ attains the minimum.
This contradicts \eqref{NMC2'}, hence $v=0$ and the claim follows.
 \end{proof}
 
\subsection{{Maximum} Principles for \eqref{QAP}}\label{Ssec:52}
  Using the formulation \eqref{NMC2} for $H_\alpha(u)$, we establish  next
  the uniqueness of solutions to \eqref{QAP} and to \eqref{eqqu} within the regularity class \eqref{reggg} as announced in the proof of Theorem~\ref{MT1}.
 \begin{lemma}\label{L:US}
 Let $\lambda_1>0$, $\lambda_2\in\R^n$, and $T>0$. If
 \begin{equation*} 
 u_1,\, u_2\in {\rm C}([0,T),{\rm h}^{1+\beta}(\bT^n))\cap {\rm C}((0,T),{\rm h}^{1+\gamma}(\bT^n))\cap {\rm C}^1((0,T),{\rm h}^{\gamma-\alpha}(\bT^n))
  \end{equation*} 
  are solutions to
   \begin{equation}\label{sps}
\frac{du}{dt}(t)=\lambda_1\Phi(u(t))[u(t)]+\lambda_2\cdot\nabla u (t),\quad t\in(0,T),\qquad  u(0)=u_0,
 \end{equation}
 then $u_1=u_2$.
 \end{lemma}
 \begin{proof}
  Set $v:=u_1-u_2$ and assume  for contradiction that $v\neq0$ in $[0,T']\times\R^n$ for some~${T'\in(0,T)}$.
Therefore, given $\e>0$, also the function $$w(t,x):=e^{-t\e}v(t,x),\quad (t,x)\in[0,T']\times\R^n,$$ is not constant.
Observing that  $w\in{\rm C}([0,T']\times \R^n)$, we may assume 
 that $w$ attains a positive maximum at ${(t_0,x_0)\in(0,T']\times\R^n}$; that is $w(t,x)\leq w(t_0,x_0)$ for all $(t,x)\in [0,T']\times \R^n.$ 
 Since $w\in {\rm C}^1((0,T)\times\R^n)$ we deduce that $\p_tw(t_0,x_0)\geq0$ and~${\nabla w(t_0,x_0)=\nabla v(t_0,x_0)=0}$.
These relations and the fact that $u_1$ and $u_2$ are solutions to \eqref{sps} lead to
\begin{align*}
0\leq \p_tw(t_0,x_0) =-\e w(t_0,x_0)+e^{-\e t_0}\p_tv(t_0,x_0)<e^{-\e t_0}\p_tv(t_0,x_0),
\end{align*}
where 
\begin{align*}
\p_tv(t_0,x_0)&=\lambda_1(1+|\nabla u_1|^2)^{1/2}(t_0,x_0)\big(H_\alpha(u_2(t_0))(x_0)-H_\alpha(u_1(t_0))(x_0)\big),
\end{align*}
and, in view of \eqref{NMC2},
\begin{align*} 
&H_\alpha(u_2(t_0))(x_0)-H_\alpha(u_1(t_0))(x_0)\\[1ex]
&=\PV\int_{\R^n}  \frac{\delta_{[x_0,y]}v(t_0)}{|y|^{n+1+\alpha}}\int_0^1  F'\Big(  (1-s)\frac{\delta_{[x_0,y]}u_2(t_0)}{|y|}+s\frac{\delta_{[x_0,y]}u_1(t_0)}{|y|} \Big)\, ds \, dy\\[1ex]
&\quad +\PV\int_{\R^n} \frac{\delta_{[x_0,y]}v(t_0)}{|y|^{n+1+\alpha}}\int_0^1  F'\Big(  -(1-s)\frac{\delta_{[x_0,y]}u_2(t_0)}{|y|}-s\frac{\delta_{[x_0,y]}u_1(t_0)}{|y|} \Big)\, ds \, dy.
 \end{align*}
 Since $F'<0$ and $v(t_0,x_0)\geq v(t_0,x_0-y)$ for all $y\in\R^n$, the latter expression is non-positive, and we obtain a contradiction.
 \end{proof}
 }
 
 {In Lemma~\ref{L:9}  we present  several maximum principles for  the evolution problem~\eqref{QAP}.}
A similar result  {can be found}, e.g., in~\cite[Lemma 5.3]{AFW22xx},
 but since these features are very special properties of the nonlocal mean curvature flow we recall them here.

\begin{lemma}\label{L:9} Let $u_0\in{\rm h}^{1+\beta}(\bT^n)$. Then, the maximal solution $u:[0,T^+)\to {\rm h}^{1+\beta}(\bT^n)$  to~\eqref{QAP} found in Theorem~\ref{MT1} satisfies
\begin{itemize}
\item[(i)] $\|u(t)\|_0\leq \|u_0\|_0$ for $ t\in[0,T^+)$.
\item[(ii)] $\|\p_{x_i}u(t)\|_0\leq \|\p_{x_i}u_0\|_0 $ for $t\in[0,T^+)$ and $1\leq i\leq n$. 
\item[(iii)]  {If $\beta>\alpha,$ then $\|\p_t u(t)\|_0\leq \| \p_tu (0)\|_0$ for $ t\in[0,T^+)$}. 
\end{itemize}
 \end{lemma}
 \begin{proof}
 Since the proof of (ii) and (iii) follows by differentiating \eqref{QAP} (with $H_\alpha(u)$ as expressed in \eqref{NMC2})
  with respect to $x_i$ and $t$, respectively, and by arguing as in the proof of (i), we only present the proof of (i) in detail. 
  
Concerning (i), we  {show} that  $\|u(t)\|_0 {<\|u_0\|_0+\e}$ for all $t\in[0,T^+)$ and $\e\in(0,1)$. 
 To prove this claim we argue by  contradiction and assume there  was $\e\in(0,1)$ such that 
  \[
  \inf\{t\in[0,T^+)\,:\|u(t)\|_0\geq \|u_0\|_0+\e \}=:t_0\in(0,T^+).
  \]
  Then $\|u(t_0)\|_0= \|u_0\|_0+\e$, hence there exists $x_0\in\R^n$ such that $|u(t_0,x_0)|= \|u_0\|_0+\e$.
  Without loss of generality we may assume that $u(t_0,x_0)>0,$ hence $u$ attains the global maximum in $[0,t_0]\times\R^n$  at $(t_0,x_0).$
Taking   $\p_tu(t_0,x_0)\geq0$ and $\nabla u(t_0,x_0)=0$  into account we get 
\[
0\leq \p_tu(t_0,x_0)=\PV\int_{\R^n} \frac{1}{|y|^{n+\alpha}}\bigg[ F \Big(  \frac{\delta_{[x_0,y]}u(t_0)}{|y|} \Big)-F \Big( - \frac{\delta_{[x_0,y]}u(t_0)}{|y|} \Big)\bigg]\, dy<0,
\]
provided that $u(t_0)$ is not constant.
The strict inequality in the latter relation   follows (when~$u(t_0)$ is not constant) by arguing as Lemma~\ref{L:8}.
Hence, if $u(t_0)$ is not constant we obtain a contradiction.

If $u(t_0)$ is constant, then there exists
\[
  t_1:=\inf\{t\in[0,T^+)\,:\|u(t)\|_0\geq \|u_0\|_0+\e/2 \}\in(0, t_0).
  \]
On the one hand, if $u(t_1)$ is not constant, we may argue as above to obtain a contradiction. On the other hand, if $u(t_1)$ is constant, then 
$u(t)=u(t_1)$ for all $t\in (t_1, T^+)$ as constants are stationary solutions,
 hence also~${u(t_0)=u(t_1)}$, which is again a contradiction.
This proves the claim.  
 \end{proof}
 
 \subsection{Proof of Theorem~\ref{MT2}}\label{Ssec:53'}
We first introduce some notation. Given 
an integrable function $u:\bT^n\to\R,$ we denote by $\langle u\rangle $ the integral mean of $u$; that is,
\[
\langle u\rangle:=\frac{1}{(2\pi)^n} \int_{\bT^n} u\, dx.
\]
Furthermore, we let $P$  be the  projection 
\[
Pu:=\langle u\rangle.
\]
Then, given $r\geq0$, we can represent ${\rm h}^r(\bT^n)$ as a direct sum
\[
{\rm h}^r(\bT^n)=P{\rm h}^r(\bT^n)\oplus(I-P){\rm h}^r(\bT^n),
\]  
where $(I-P){\rm h}^r(\bT^n)$ is the subspace of ${\rm h}^r(\bT^n)$   consisting of those functions with zero integral mean.

\begin{proof}[Proof of Theorem~\ref{MT2}]
Let $u:[0,T^+)\to{\rm h}^{1+\beta}(\bT^n)$ be the maximal solution to 
 \begin{equation*} 
 \frac{du}{dt}(t)=\Phi(u(t))[u(t)],\qquad t>0, 
 \end{equation*}
 determined by $ u(0)=u_0\in {\rm h}^{1+\beta}(\bT^n)$,  as found in Theorem~\ref{MT1}.
 Given $t\in[0,T^+)$, let
\begin{equation}\label{NEVS}
q(t):=\langle u(t)\rangle \qquad\text{and}\qquad  v(t):=u(t)-q(t).
 \end{equation}
 Then $v(0)=u_0-\langle u_0\rangle \in(I-P){\rm h}^{1+\beta}(\bT^n)$  and $v(t)$ and $u(t)$ differ just by a $t$-dependent constant if $t>0$.
 Recalling the definition of $\Phi$, we therefore have~$\Phi(u(t))[u(t)]=\Phi(v(t))[v(t)]$ for~$t\in[0,T^+)$.
 The latter observation, the properties of a solution to \eqref{QAP}, and \eqref{NEVS} lead us now to
 \begin{align*}
 \frac{dv}{dt}(t)&= \frac{du}{dt}(t)-\Big\langle  \frac{du}{dt}(t)\Big\rangle=\Phi(u(t))[u(t)]-\langle \Phi(u(t))[u(t)]\rangle\\[1ex]
 &=\Phi(v(t))[v(t)]-\langle \Phi(v(t))[v(t)]\rangle,\quad t>0.
 \end{align*}
We now define the operator
\begin{equation*}
\Psi(u)[v]:=\Phi(u)[v]-\langle \Phi(u )[v]\rangle,
\end{equation*}
and infer from \eqref{MP} that
\begin{equation*} 
\Psi\in{\rm C}^\infty\big((I-P){\rm h}^{1+\beta}(\bT^n),\mathcal{L}(E_1, E_0)\big),
\end{equation*}
where
 {
\[ 
E_0:=(I-P){\rm h}^{\gamma-\alpha}(\bT^n)\qquad \text{and}\qquad E_1:=(I-P){\rm h}^{1+\gamma}(\bT^n).
 \] }
With this notation, we have identified  $v$ as a solution to the  evolution problem
 \begin{equation}\label{QAP'S}
\frac{dv}{dt}(t)=\Psi(v(t))[v(t)], \quad t>0,\qquad v(0)=v_0:=u_0-\langle u_0\rangle,
\end{equation}
 which is quasilinear and of parabolic type. Indeed, with respect to the decompositions 
 \begin{align*}
 {\rm h}^{1+\gamma}(\bT^n)=P{\rm h}^{1+\gamma}(\bT^n)\oplus  {E_1} \qquad \text{and}\qquad
{\rm h}^{\gamma-\alpha}(\bT^n)=P{\rm h}^{\gamma-\alpha}(\bT^n)\oplus {E_0},
 \end{align*}
we can represent $\Phi(u)$ as a matrix operator
\[
\Phi(u)=
\begin{pmatrix}
0&\langle \Phi(u)[\cdot]\rangle\\[1ex]
0&\Psi(u)
\end{pmatrix}.
\] 
This representation  together with Theorem~\ref{T:GP} and   \cite[Corollary I.1.6.3]{Am95} implies that 
\begin{equation*}
-\Psi\in{\rm C}^\infty\big((I-P){\rm h}^{1+\beta}(\bT^n),\mathcal{H}(E_1, E_0)\big).
\end{equation*}
Moreover, $v_*:=0$ is a stationary solution to \eqref{QAP'S} and the linearized operator
\[
\bA:=\Psi(v_*)+(\p\Psi(v_*)[\cdot])[v_*]=\Psi(v_*)\in\mathcal{L}(E_1,E_0)
\] 
is, in view of \eqref{phi0}, \eqref{fmult}, \eqref{symbol}, and \eqref{repr} the Fourier multiplier with symbol 
\[
m(k)=-\omega_0|k|^{1+\alpha}, \qquad k\in\mathbb{Z}^n\setminus\{0\},
\]
where the negative constant $-\omega_0$ is defined by the integral \eqref{repr}  with $a=0$ (which is  thus constant in $x$).
Since the embedding $E_1 \hookrightarrow E_0$ is compact, 
we infer from \cite[Theorem~III.6.29]{Ka95} that the spectrum $\sigma(\bA)$ consists entirely of isolated eigenvalues with finite algebraic multiplicities.
 Consequently,
\[
\sigma(\bA)=\{-\omega_0|k|^{1+\alpha}\,:\, k\in\mathbb{Z}^n\setminus\{0\}\}.
\]
Hence, $\re\lambda\leq -\omega_0<0$ for all $\lambda\in\sigma(\bA)$ and we are in a position to apply the quasilinear principle of linearized stability  \cite[Theorem 1.3]{MW20} in the context of \eqref{QAP'S}.
Therefore, given~${\omega\in(0,\omega_0)}$, there exist constants $\e>0$ and $M\geq1$ such that for each $\|v_0\|_{1+\beta}\leq \e$, the solution $v=v(\cdot; v_0)  $ to \eqref{QAP'S} exists globally and 
\begin{equation}\label{stestS}
\|v(t)\|_{1+\beta}\leq M e^{-\omega t}\|v_0\|_{1+\beta},\qquad t\geq0.
\end{equation}
As already mentioned in the discussion subsequent to Theorem \ref{MT1}, we point out again that
 the nonlocal mean curvature of order $\alpha$ is not well-defined for $\beta\leq \alpha$. 
 The property \eqref{stestS} can thus not be used directly to derive estimates for~${\Phi(v(t))[v(t)]}$, $t>0$. Therefore,
in order to establish the convergence of $q(t)$ towards a constant we revisit the proof of \cite[Theorem~1.3]{MW20} to estimate $\|v(t)\|_{1+\gamma}$ for $t>0$.

To start with, we set
\[
A(t):=-\Psi(v(t)),\qquad t\geq0.  
\]
Similarly as in the proof of Theorem~\ref{MT1}, we choose a H\"older exponent $\beta'\in(0,\beta)$  such that~$\gamma<\alpha+\beta'$.
We further set $\rho:=(1+\alpha)^{-1}(\beta-\beta')$ and  $4\delta:=\omega_0-\omega$.
As shown in the proof of \cite[Theorem 1.3]{MW20}, if $\|v_0\|_{1+\beta}\leq \e$, the solution  $v=v(\cdot; v_0)  $ to~\eqref{QAP'S} additionally satisfies
\begin{equation}\label{exttex}
\|A(t)-A(s)\|_{\kL(E_1,E_0)}\leq c\|v_0\|_{1+\beta}|t-s|^\rho,\qquad\text{$t,\,s\in[0,\infty),$}
\end{equation}
and
\begin{equation}\label{exttex2}
-\omega_0+\delta+A(t)\in\mathcal{H}(E_1, E_0,\kappa,\delta),\qquad\text{$t\in[0,\infty),$}
\end{equation} 
 with a fixed constant $\kappa\geq1$, see Section~\ref{Sec:1.1}. 

In view of \eqref{exttex}-\eqref{exttex2} and the interpolation relation~\eqref{interpol} we may apply \cite[Lemma~II.5.1.3]{Am95} to deduce that the  evolution operator~$U_A$ associated  with the mapping 
$$[t\mapsto A(t)]\in {\rm C}^\rho([0,\infty), \kL(E_1,E_0))$$ 
 satisfies, after choosing a smaller $\e$ if necessary, the estimate
\begin{equation*}
\|U_A(t,s)\|_{\kL((I-P){\rm h}^{1+\beta}(\bT^n), E_1)}\leq ce^{-\omega(t-s)}(t-s)^{-\rho'},\qquad 0\leq s<t<\infty,
\end{equation*}
where $\rho':=(1+\alpha)^{-1}(\gamma-\beta).$ 
From this and \cite[Remark II.2.1.2]{Am95} we deduce that 
\begin{equation}\label{hgh}
\|v(t)\|_{1+\gamma}=\|U_{A}(t,0)v_0\|_{1+\gamma}\leq  ce^{-\omega t}t^{-\rho'}\|v_0\|_{1+\beta},\qquad t>0. 
\end{equation}
We now consider  the function $q$ representing the integral mean of $u$ for which we have
\[
q'(t)=\Big\langle  \frac{du}{dt}(t)\Big\rangle=\langle \Phi(u(t))[u(t)]\rangle=\langle \Phi(v(t))[v(t)]\rangle,\quad t>0.
\]
This implies that  $q(t)$ exists for all  $t\geq0$.
Since $q$ is continuous on $[0,\infty)$ we also have
\[
\langle u_0\rangle=q(0)=\lim_{t\to 0}q(t)
\]
 and hence
\begin{align}\label{eq0}
\langle u(t)\rangle=q(t)=\langle u_0\rangle+\int_0^t\langle\Phi(v(s))[v(s)]\rangle\, ds, \qquad t\geq0.
\end{align}
Indeed,  the definition of $\Phi$, \eqref{f1}, and \eqref{hgh} imply for $x\in\R^n$  and $t> 0$ that
\begin{align*}
|\Phi(v(t))[v(t)](x)|&\leq  c\int_{\R^n} \frac{|\delta_{[x,y]}v(t)-y\cdot\nabla v(t,x-y)|}{|y|^{n+1+\alpha}}dy\leq ce^{-\omega t}t^{-\rho'}\|v_0\|_{1+\beta},
\end{align*}
and therefore
\begin{equation}\label{eq1}
\int_0^\infty|\langle\Phi(v(s))[v(s)]\rangle|\, ds\leq c\|v_0\|_{1+\beta}\int_0^\infty e^{- \omega s}s^{-\rho'}\, ds \leq c\|v_0\|_{1+\beta}.
\end{equation}
This ensures in particular the convergence of the integral in \eqref{eq0}.
Moreover, given $t\geq 1,$ we may estimate 
\begin{equation}\label{eq2}
\int_t^\infty|\langle\Phi(v(s))[v(s)]\rangle|\, ds\leq c\|v_0\|_{1+\beta}\int_t^\infty e^{- \omega s} \, ds \leq ce^{-\omega t}\|v_0\|_{1+\beta}.
\end{equation}
 Introducing the constant  
\begin{equation}\label{cu0}
C(u_0):=\lim_{t\to\infty}q(t)=\langle u_0\rangle+\int_0^\infty\langle\Phi(v(s))[v(s)]\rangle\, ds 
\end{equation}
 we infer  from \eqref{eq0}-\eqref{eq2} that 
\begin{equation}\label{eq7}
|\langle u(t)\rangle-C(u_0)|\leq \int_t^\infty |\langle\Phi(v(s))[v(s)]\rangle|\, ds\leq M e^{-\omega t}\|u_0-\langle u_0\rangle\|_{1+\beta},\qquad t\geq0,
\end{equation}
for a possibly larger constant $M\geq1$. 
Gathering \eqref{NEVS}, \eqref{stestS}, and \eqref{eq7} implies \eqref{Exst}.
 Finally, Lemma~\ref{L:9}~(i) ensures that $|C(u_0)|\le \|u_0\|_0$. 
 This yields Theorem~\ref{MT2}.
\end{proof}

\appendix
\section{A Lifiting Property for Periodic Besov Spaces on $\bT^n$}\label{Sec:A}
 We recall the definition of Besov spaces and refer to \cite{ST87} for more details. As usual, $ \mathcal{D}'(\bT^n)$ is the topological dual of~$\mathcal{D}(\bT^n)={\rm C}^\infty(\bT^n)$, and any $f\in  \mathcal{D}'(\bT^n)$ can be represented as
\[
f=\sum_{k\in\Z^n} \wh f(k)e^{ik\cdot x},
\]
where
\[
\wh f(k):=(2\pi)^{-n}f(e^{-ik\cdot x}),\qquad k\in\Z^n.
\]
Let $(\varphi_j)_{j\in\N}\subset \mathcal{S}(\R^n)$ be such that 
$$
\supp\varphi_0\subset\{|x|\leq 2\}\,,\quad  \supp\varphi_j\subset\{2^{j-1}\leq |x|\leq 2^{j+1}\}
$$ 
for $j\geq1$ with
\[
\sum_{j=0}^\infty \varphi_j=1\ \text{ on }\ \R^n
\]
and, for each $\mu\in\N^n$, 
\[
\sup_{j\geq0}\sup_{x\in\R^n}2^{j|\mu|}|\p^\mu \varphi_j(x)|<\infty.
\]
Given $s\in\R$ and $1\le p,q\le \infty$, the Besov space $B^{s}_{p,q}(\bT^n)$ consists of all $f\in  \mathcal{D}'(\bT^n)$ for which the norm 
\begin{equation}\label{bspq}
\|f\|_{B^{s}_{p,q}(\bT^n)}:=\left\|\left( 2^{sj}\bigg\|\sum_{k\in\Z^n}\varphi_j(k) \wh f(k) e^{ik\cdot x}\bigg\|_{L_p(\bT^n,dx)}\right)_{j\in\N}\right\|_{\ell_q}
\end{equation}
is finite. Then $B^{s}_{p,q}(\bT^n)$ equipped with the norm $\|\cdot\|_{B^{s}_{p,q}(\bT^n)}$ is a Banach space (with equivalent norms when replacing $(\varphi_j)_{j\in\N}\subset \mathcal{S}(\R^n)$ by another family with the same properties). For our purposes it is important to note that one may identify ${\rm C}^s(\bT^n)$
with the Besov space~${B^{s}_{\infty,\infty}(\bT^n)}$ provided that~${s>0}$ with $s\not\in\N$.

As a preparatory result we provide in Lemma~\ref{L:A1} a direct proof for
the boundedness of the Fourier multiplier $I_t$ defined by the symbol $\mathfrak{m}_t(k)=|k|^{t} $ for $k\neq0 $ and $\mathfrak{m}_t(0)=1$ from $B^{s}_{p,q}(\bT^n)$ to 
$B^{s-t}_{p,q}(\bT^n)$ when $s,\, t\in\R$. 
In particular, Lemma~\ref{L:A1} shows that $I_t$ is an isomorphism with inverse $I_{-t}$. A similar lifting property is established in the
 nonperiodic case in \cite[Theorem 2.3.8]{Tr83}. As we are lacking a precise reference for the periodic case, we include a proof here.

\begin{lemma}\label{L:A1} If $s,\, t\in\R$ and  $1\le p,q\le \infty$, then
\[
I_t:=\bigg[\sum_{k\in\Z^n}  \wh f(k) e^{ik\cdot x}\mapsto \sum_{k\in\Z^n}\mathfrak{m}_t(k) \wh f(k) e^{ik\cdot x}\bigg]\in\kL(B^{s}_{p,q}(\bT^n),B^{s-t}_{p,q}(\bT^n)),
\]
where $\mathfrak{m}_t(k)=|k|^{t} $ for $k\neq0 $ and $\mathfrak{m}_t(0)=1$.
\end{lemma}

 The proof of Lemma~\ref{L:A1} is very much inspired by a result established for operator valued Fourier multipliers in the case $n=1$, see \cite{AB04}.
 Following \cite[Proposition 2.2]{AB04}, we prove the following auxiliary result which is the basis for the proof of Lemma~\ref{L:A1}.

\begin{lemma} \label{L:A2}  If $\zeta\in \mathcal{S}(\R^n)$ and $1\le p\le\infty$,  then
\begin{equation*}
\bigg\|\sum_{|k|\leq M}  \zeta (k) \wh f(k) e^{ik\cdot x}\bigg\|_{ L_p(\bT^n)}\leq (2\pi)^{-n/2}\|\mathcal{F}^{-1}\zeta\|_{L_1(\R^n)} \|f\|_{L_p(\bT^n)}
\end{equation*}
for every trigonometric polynomial $f=\sum_{|k|\leq M}\wh f(k) e^{ik\cdot x}$ with $M\in\N$.
\end{lemma}

\begin{proof}
Since the (nonperiodic) Fourier transform $\mathcal{F}:\mathcal{S}(\R^n)\to \mathcal{S}(\R^n)$ is a bijection, we have
\begin{align*}
\bigg\|\sum_{|k|\leq M}\zeta(k) \wh f(k) e^{ik\cdot x}\bigg\|_{L_p(\bT^n)}&=(2\pi)^{-n/2}\bigg\| \int_{\R^n} \mathcal{F}^{-1}\zeta(\xi) \sum_{|k|\leq M}\wh f(k)   e^{ik\cdot (x-\xi)}\, d\xi\bigg\|_{L_p(\bT^n)}\\[1ex]
&=(2\pi)^{-n/2}\bigg\| \int_{\R^n} \mathcal{F}^{-1}\zeta(\xi) f(x-\xi)\, d\xi\bigg\|_{L_p(\bT^n)},
\end{align*}
and the claim follows from Young's inequality.
\end{proof}

We are now in a position to establish Lemma~\ref{L:A1}.

\begin{proof}[Proof of Lemma~\ref{L:A1}]
We may assume that $\supp\varphi_0\subset\{|x|\leq \sqrt{2}\}$. Noticing that~${\mathfrak{m}_t(k)=1}$ if $|k|\leq \sqrt2$,  we have for $j=0$
\begin{align*}
2^{(s-t)j}\bigg\|\sum_{k\in\Z^n}\varphi_0(k) \mathfrak{m}_t(|k|)\wh f(k) e^{ik\cdot x}\bigg\|_{L_p(\bT^n)}&= 2^{sj}\bigg\|\sum_{|k|\leq 1}\varphi_0(k) \wh f(k) e^{ik\cdot x}\bigg\|_{L_p(\bT^n)}.
\end{align*}

We then estimate the terms with $j\geq 1$. 
To this end we choose a function $\psi\in\mathcal{S}(\R)$ such that $\supp\psi\subset\{1/4\leq |x|\leq4\}$ and 
$\psi=1$ on $\{1/2\leq |x|\leq2\}$. 
Noticing that  the function~${\zeta:\R^n\to\R}$ with~$\zeta(x):=|x|^t\psi(|x|)$ belongs to $\mathcal{S}(\R^n)$,   Lemma~\ref{L:A2} yields
\begin{align*}
&2^{(s-t)j}\bigg\|\sum_{k\in\Z^n}\varphi_j(k) \mathfrak{m}_t(k)\wh f(k) e^{ik\cdot x}\bigg\|_{L_p(\bT^n)}\\[1ex]
&= 2^{sj}\bigg\|\sum_{2^{j-1}\leq |k|\leq 2^{j+1}}\psi(2^{-j}|k|)(2^{-j}|k|)^{t}\varphi_j(k)  \wh f(k) e^{ik\cdot x}\bigg\|_{L_p(\bT^n)}\\[1ex]
&\leq (2\pi)^{-n/2}\big\|\mathcal{F}^{-1}[x\mapsto \zeta(2^{-j}x)]\big\|_{L_1(\R^n)} 2^{sj}\bigg\|\sum_{k\in\Z^n}\varphi_j(k)  \wh f(k) e^{ik\cdot x}\bigg\|_{L_p(\bT^n)}\\[1ex]
&\leq (2\pi)^{-n/2}\big\|\mathcal{F}^{-1}\zeta\big\|_{L_1(\R^n)}\, 2^{sj}\bigg\|\sum_{k\in\Z^n}\varphi_j(k)  \wh f(k) e^{ik\cdot x}\bigg\|_{L_p(\bT^n)},
\end{align*}
where in the last step we used the  property  that $$\|\mathcal{F}^{-1}\zeta\big\|_{L_1(\R^n)}=\|\mathcal{F}^{-1}[x\mapsto \zeta(b x)]\|_{L_1(\R^n)}$$ for~$ b>0$.
Recalling~\eqref{bspq}, the assertion follows.
\end{proof}

\bibliographystyle{siam}
\bibliography{WM22}
\end{document}